# Large Deviations and Queueing Networks: Methods for Rate Function Identification*

Rami Atar and Paul Dupuis

Lefschetz Center for Dynamical Systems
Brown University
Division of Applied Mathematics
Providence, R.I. 02912

October 4, 1998

## Abstract

This paper considers the problem of rate function identification for multidimensional queueing models with feedback. A set of techniques are introduced which allow this identification when the model possesses certain structural properties. The main tools used are representation formulas for exponential integrals, weak convergence methods, and the regularity properties of associated Skorokhod Problems. Two examples are treated as special cases of the general theory: the classical Jackson network and a model for processor sharing.

## 1 Introduction

Although there has been considerable interest in establishing a theory of large deviations for queueing networks and related systems [17], there are few general results for multidimensional systems with feedback, save theorems which establish the existence of a large deviation principle but fail to provide an explicit formula for the rate function [3]. In this paper we introduce techniques that allow one to fill in this gap, at least for families of networks that possess certain structural properties. The main tools we use are the representation formulas used to prove existence in [3], weak convergence methods, and the regularity properties of an associated Skorokhod Problem.

We begin with a review of the literature that deals with large deviation properties of queueing networks. One of the first papers on the topic is [8], which establishes certain large deviation properties of Jackson networks via nonlinear PDE techniques. Unfortunately, the methods of this paper do not extend easily to more general situations. A probabilistic method is used to prove large deviation upper bounds for a fairly general class of Markov models in [6]. The corresponding lower bound is not proved, and so the tightness of these upper bounds remains an open question, although a partial answer will be given in the present work.

---

*This research was supported in part by the National Science Foundation (NSF-DMS-9704426) and the Army Research Office (DAAH04-96-1-0075).

A number of techniques have been developed that are very much tailored to particular models. For example, the special class of tandem queues is quite tractable, in large part because the absence of feedback means that continuous mapping methods can be applied. Indeed, in this case one can represent the queueing model as the composition of the mapping on path space defined by a suitable Skorokhod Problem and an unconstrained process [7]. When such a representation is available the large deviations analysis can be carried out by applying the contraction principle. Results for models of this type can be found in [2, 14, 18]. Although this method can be extended to cover a broader class of models (e.g., feedforward as in [13, 15]), it breaks down when feedback is present. A general result on large deviations for processes whose statistical behavior can be discontinuous across a smooth $(n-1)$-dimensional interface in $\mathbb{R}^n$ was proved in [5], and then applied in [11] to prove the large deviation principle for a general class of stable two dimensional queueing models. Related results that also rely on a reduction to what are essentially one dimensional problems include [1, 19].

The only existing theory that does not make significant use of model specific geometric features is presented in [3], which considers a general class of jump Markov processes that model queueing systems. This paper proves the existence of a large deviation principle, and also provides a characterization of the rate function. The paper falls short, however, in that it does not identify the rate function.

As the previous discussion suggests, no general method has emerged for rate function identification when there is feedback and if the dimension is greater than two. On the other hand, an explicit expression for an upper large deviation rate function may be found in [6]. As we will see below, for many models the upper bound proved in [6] is actually tight, and thus the main difficulty appears in the proof of the large deviation lower bound. There are several reasons for this difficulty. One has to do with the fact that queueing systems fall into the category of "processes with discontinuous statistics," as defined in [5, 6]. The discontinuities appear because the generator of a queueing process often changes abruptly when one or more of the queues becomes empty.

Before discussing in detail how these discontinuities affect the analysis, it is important to note that our proofs will actually be based on a control theoretic representation for the large deviation probabilities. This approach has much in common with the change of measure argument often used to prove large deviation lower bounds, and in fact each control will correspond to one such change of measure. Since the change of measure argument is more widely known in this context, we will use the terminology of this technique instead, with an understanding that the there is an equivalent phrasing in terms of control representations.

Thus we resume our consideration of how one may establish a large deviation lower bound, and assume that a change of measure has been selected for this purpose. For processes with discontinuous statistics it turns out that one must characterize certain properties of the asymptotic fractions of time that the process spends in each subregion of smooth statistical behavior. The situation is simplest when one can uniquely characterize the fractions of time themselves. However, this is not usually possible if the dimension is greater than two, and it is this difficulty which has thwarted a probabilistic proof of even the (relatively) simple case of Jackson networks.

A method for eliminating this difficulty is one of our main innovations. We show, for the models introduced later in the paper, that an argument based on Jensen's inequality allows one to restrict the changes of measure that one must consider. In particular, it turns out that given a network from among the classes we consider (e.g., a Jackson network), one can restrict to changes of measure that return the process to this same class. This is extremely useful when analyzing the asymptotic properties of the new process, and the reason has to do with the second main new ingredient we use: weak convergence methods and the Skorokhod Problem. For the types of models we treat, it is known that under the standard large deviation scaling (which is the same as the law of large numbers (LLN) scaling), limits of the queueing system can be characterized as the unique solution of a well behaved Skorokhod Problem [9]. Recall that the use of Jensen's inequality discussed above allows one to restrict, a priori, to changes of measure that define processes of the same sort as the original process. Because they fall into the same class as the original model, the LLN limit of each model defined by one of these changes of measure will be characterized by a well behaved SP (which is typically *not* the SP associated with the original model). This is rather convenient, since the uniqueness of solutions to a well behaved SP is just what is needed to properly characterize the limiting behavior of the process.

Besides a well behaved Skorokhod Problem, an additional structural property that is needed is related to the form of the rate function. As noted previously, the formulation of the rate function involves quantities that may be interpreted as limits of fractions of time spent in different subregions of smooth statistical behavior. In general, the limits of the fractions of time are not uniquely characterized, since there are too few constraints when compared to the number of subregions of smooth behavior. However, in the course of applying Jensen's inequality we derive an alternative representation for the rate function. This representation involves only functionals of the fractions of time that are uniquely characterized, and thus weak convergence methods can be applied to prove the desired convergence.

As an applications of these methods we identify the rate function for two families of continuous time jump Markov models. However, the range of applications is somewhat broader, and for example one can consider also Markov modulated jump rates, discrete time processes, and other variations. The most significant restriction appears to be the requirement that the associated SPs must all be regular.

An outline of the paper is as follows. In Section 2 notation is introduced and the needed background material is presented. It is convenient to introduce "localized" versions of the original queueing model. These models provide, for a given point in the state space of the original model, the simplest model whose large deviations behavior is the same as that of the original model near the given point. We then state the control representation for the rate function for the localized models, and introduce the Skorokhod Problem, which will be used to identify weak limits in the asymptotic analysis of the representation formula. In Section 3 we relate various representations for the rate function for the local model, which will eventually appear as the integrand in the rate function for the full model. It is in this section that we isolate key parameters in the rate function that are needed for its identification. In Section 4 a number of abstract assumptions are made regarding these parameters and other properties of the network, and weak convergence methods are applied to complete the identification of the rate function. Finally, in Section 5 we verify the assumptions for two families of models. Some minor details are relegated to an Appendix to ease the exposition.

## 2 Background

### 2.1 Stochastic control representation and statement of the LDP

This subsection summarizes some definitions and results of [3]. Under appropriate assumptions, the large deviation principle holds for a general class of jump Markov processes, and the rate function may be characterized in terms of rate functions for local models, as explained below (Theorems 2, 3). The proof of the large deviation principle is based on a control-theoretic representation for the probabilities that these processes stay within tubes centered at piecewise linear paths (Theorem 1).

For $T = [0, t]$ or $T = [0, \infty)$ and a Polish space $\mathcal{S}$, we denote by $D(T : \mathcal{S})$ the space of all càdlàg functions $T \mapsto \mathcal{S}$. The process used to model the queueing system is a jump Markov process on the orthant $\mathbb{Z}_+^N = \{x \in \mathbb{Z}^N : \langle x, e_i \rangle \geq 0, i = 1, 2, \ldots, N\}$, with paths in $D([0, \infty) : \mathbb{Z}_+^N)$. (In [3] more general polyhedra than the orthant are allowed.) The analysis in [3] relies on the introduction of local models and the stochastic processes associated with them. These processes take values in subsets of $\mathbb{Z}^N$ of the form $\mathbb{Z}_+^{N,K}$, where $K \subset \{1, 2, \ldots, N\}$, and

$$\mathbb{Z}_+^{N,K} \doteq \{x \in \mathbb{Z}^N : \langle x, e_i \rangle \geq 0, i \in K\}.$$

Thus for a local model associated with the state space $\mathbb{Z}_+^{N,K}$, non-negativity constraints are enforced on $x_i, i \in K$, but not otherwise. The original queueing model is associated with $\mathbb{Z}_+^{N,\{1,\ldots,N\}} = \mathbb{Z}_+^N$. The term *full model* will be used to distinguish the original queueing model from its related family of local models. We denote also

$$\mathbb{R}_+^N \doteq \{x \in \mathbb{R}^N : \langle x, e_i \rangle \geq 0, i = 1, 2, \ldots, N\},$$

$$\mathbb{R}_+^{N,K} \doteq \{x \in \mathbb{R}^N : \langle x, e_i \rangle \geq 0, i \in K\}.$$

As will be discussed further below, the large deviation behavior of a queueing model can often be determined by considering the asymptotics of the probability that the process stays in a small neighborhood of a constant velocity trajectory, if we allow the initial position and velocity of the trajectory to range over suitable values. Suppose, for example, that the trajectory is $y + \beta t$, $t \in [0, 1]$. If the initial point $y$ or terminal point $y + \beta$ of the trajectory satisfy more active constraints of the form $\langle e_i, x \rangle = 0$ than the "interior" of the trajectory $\{x = y + \beta t : t \in (0, 1)\}$, then the more complicated dynamics in the neighborhood of such an end point are an annoying but relatively unimportant nuisance. See Figure 1.

The appropriate local model will effectively throw away all parts of the process that are unimportant when determining large deviation properties of the full model with regard to such constant velocity trajectories. For more details, see the discussion in [3, Section 4] and Example 1 below.

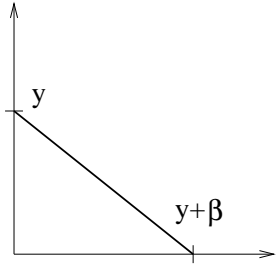

Figure 1: A trajectory whose endpoints satisfy more constraints than the interior

The notation introduced in a moment will be used in the context of both the original and the local models.

We first introduce the notion of a facet. Facets will be the regions of constant statistical behavior, and they will vary with the local model under consideration. For $K \subset \{1, 2, \ldots, N\}$ and $I \subset K$, let
$$\mathcal{F}^{K,I} \doteq \{x \in \mathbb{R}^N : \langle x, e_j \rangle = 0, j \in I, \langle x, e_i \rangle > 0, i \in K \setminus I\}.$$
Recall that for a local model associated with $K$, only the variables $x_i, i \in K$ are constrained to be non-negative. Facets are characterized by the subset of indices $I$ for which the inequality constraint is tight. For each local model it will be necessary to group those states for which the generator of the process takes the same form. Each such group will coincide with the intersection of $\mathbb{Z}^N$ with one of the facets associated with the given local model. In this context, notice that for each subset $K \subset \{1, 2, \ldots, N\}$, the sets $\mathcal{F}^{K,I}$ partition $\mathbb{R}_+^{N,K}$ as $I$ ranges over all choices of $I \subset K$. Note also that every facet is a cone, since $x \in \mathcal{F}^{K,I}$ implies $\alpha x \in \mathcal{F}^{K,I}$ for all $\alpha > 0$. See Example 1 below for the details of a particular case.

For a Markov process with state space $\mathcal{S} \subset \mathbb{Z}^N$, we let $r(x, v) \geq 0$ denote the jump intensity from $x \in \mathcal{S}$ to $x + v \in \mathcal{S}$. We extend $r(\cdot, \cdot)$ to $\mathcal{S} \times \mathbb{Z}^N$ by letting $r(x, v) = 0$ for $x \in \mathcal{S}$, $x + v \in \mathbb{Z}^N \setminus \mathcal{S}$.

The following condition on a state space $\mathcal{S}$ and on an intensity function $r$ will be needed for each local model in order to obtain large deviation properties of the full model. Thus when the condition is assumed $\mathcal{S}$ and $r$ will vary with the particular local model. However, when parts 2 and 3 of the condition hold for the full model, they automatically hold for all local models as well.

**Condition 1** *There exists $K \subset \{1, 2, \ldots, N\}$ such that the following holds.*

1. $\mathcal{S} = \mathbb{Z}_+^{N,K}$,

2. *For all $I \subset K$, $r(x, v)$ is independent of $x \in \mathcal{F}^{K,I} \cap \mathcal{S}$,*

3. *For all $I \subset K$ and for every $x \in \mathcal{F}^{K,I} \cap \mathcal{S}$, the set $\{v \in \mathbb{Z}^N : r(x, v) > 0\}$ is finite.*

A consequence of Condition 1 is that $r(x, v)$ is radially homogeneous in $x$, and it is uniformly bounded above and below by positive constants on the set $\{(x, v) : r(x, v) > 0\}$.

For a fixed $K \subset \{1, 2, \ldots, N\}$ (or equivalently, for a fixed state space $\mathcal{S} = \mathbb{Z}_+^{N,K}$), the set of all facets $\mathcal{F}^{K,I}$ are indexed by all subsets $I \subset K$. The set of all "possible" jump directions for the full model (i.e., all vectors $v$ such that the jump rate from $x$ to $x+v$ is strictly positive) will be denoted by $V$. This set will automatically include all possible jump directions of all associated local models.

Let $\{X(t), t \in [0, \infty)\}$ be a jump Markov process on $\mathcal{S} \subset \mathbb{Z}^N$, with paths in $D([0, \infty) : \mathcal{S})$ and jump intensities $r(\cdot, \cdot)$, such that $\mathcal{S}$ and $r$ satisfy Condition 1. For $n \in \mathbb{N}$ let the processes $\{X^n(t), t \in [0, 1]\}$ be defined by

$$X^n(t) \doteq \frac{1}{n} X(nt),$$

and let the corresponding state spaces be denoted by

$$\mathcal{S}^n \doteq \left\{ \frac{1}{n} x : x \in \mathcal{S} \right\}.$$

For $\beta \in \mathbb{R}^N, n \in \mathbb{N}, y \in \mathcal{S}^n$, and $\epsilon \in (0, 1)$, let

$$p^n(y; \beta, \epsilon) \doteq P_y^n \left\{ \sup_{t \in [0,1]} \|X^n(t) - t\beta\| < \epsilon \right\},$$

where $P_y^n$ denotes probability conditioned on $X^n(0) = y$, and $\|\cdot\|$ denotes the Euclidean norm on $\mathbb{R}^N$. Let also

$$q^n(y; \beta, \epsilon) \doteq -\frac{1}{n} \log p^n(y; \beta, \epsilon),$$

where $-\log 0 \doteq \infty$.

All of the definitions given above are central to the localization procedure used in [3]. A rough explanation of how they are used is the following. In process level large deviations, one can often deduce the full LDP (cf. Theorem 3 below) if one knows the asymptotic behavior of the probability of staying in a small neighborhood of a given trajectory. An approximation argument can be used to show that it suffices to consider only trajectories that are piecewise linear, and one can then use the Markov property to simplify even further, and restrict attention to the time intervals on which the trajectory has constant velocity. Suppose that the time interval of interest is $[a, b]$. Now it turns out that the large deviations behavior over a short time interval, conditioned on starting at a point $y$, will depend only on the form of the generator near $y$ (see [3]). This implies that not all aspects of the full model are relevant when determining the probability that the process stays near the trajectory during such a time interval. In fact, if the relative interior of the given segment of the trajectory lies in the facet $\mathcal{F}^{\{1,\ldots,N\},I}$ of the full model, then the asymptotics of these probabilities are entirely determined by the behavior of a local model with state space $\mathbb{Z}_+^{N,I}$. If a given segment of a trajectory has velocity $\beta$ and if $X$ (respectively, $X^n$) is the appropriate local model (respectively, scaled local model), then asymptotic properties of $p^n$ identify the large deviation properties of the full model over the given time interval (cf. Theorems 2 and 3 below).

**Example 1** Consider a 3-dimensional model and a trajectory of the form $\{x = y + \beta t, t \in (0, 1)\}$. If $y_1 = y_3 = 0$, $\beta_1 = \beta_3 = 0$ and $y_2 \wedge (y_2 + \beta_2) > 0$, then the proper local model corresponds to $K = \{1, 3\}$. See Figure 2.

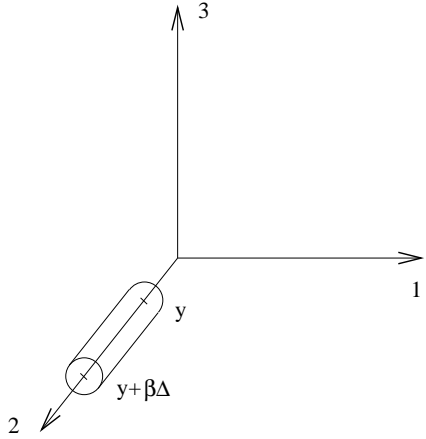

Figure 2: Localization for a three dimensional model

The four regions of the state space that are relevant are the intersections of $\mathbb{Z}^3$ with

$$\begin{aligned} A_1 &= \{(x_1, x_2, x_3) : x_1 = 0, x_2 > 0, x_3 = 0\}, \\ A_2 &= \{(x_1, x_2, x_3) : x_1 = 0, x_2 > 0, x_3 > 0\}, \\ A_3 &= \{(x_1, x_2, x_3) : x_1 > 0, x_2 > 0, x_3 = 0\}, \\ A_4 &= \{(x_1, x_2, x_3) : x_1 > 0, x_2 > 0, x_3 > 0\}. \end{aligned}$$

The trajectory is a subset of $A_1$, but small neighborhoods of the trajectory also intersect the other three sets. The sets $A_1, ..., A_4$ make up just 4 of the 8 facets of the full 3 dimensional model. All other facets are unimportant since they are a positive distance from the trajectory. The localized model has state space $\mathbb{Z}_+^{3,\{1,3\}}$, and facets

$$\begin{aligned} B_1 &= \{(x_1, x_2, x_3) : x_1 = 0, x_3 = 0\}, \\ B_2 &= \{(x_1, x_2, x_3) : x_1 = 0, x_3 > 0\}, \\ B_3 &= \{(x_1, x_2, x_3) : x_1 > 0, x_3 = 0\}, \\ B_4 &= \{(x_1, x_2, x_3) : x_1 > 0, x_3 > 0\}. \end{aligned}$$

To analyze the probability that a local model follows a linear trajectory, one can consider a representation for the pre-limit objects $q^n(y; \beta, \epsilon)$ in terms of a stochastic control problem (see also [4] for a full exposition of this approach). We recall this representation below in Theorem 1. To state it we need a few definitions.

Let $n \in \mathbb{N}$ be fixed. A *control* $u^n(y, v, t)$ is a measurable function mapping $\mathcal{S}^n \times \mathbb{Z}^N \times [0, 1]$ into $[0, \infty)$. We will impose the following condition on a control: for all $(y, v, t) \in \mathcal{S}^n \times \mathbb{Z}^N \times [0, 1]$

$$r(y, v) = 0 \text{ implies } u^n(y, v, t) = 0. \tag{1}$$

A Markov process $\{\Xi^n(t), t \in [0, 1]\}$ with state space $\mathcal{S}^n$ and jump intensity $u^n(y, v, t)$ from $y$ to $y + v/n$ at time $t$ (with $y \in \mathcal{S}^n$, $v \in \mathbb{Z}^N$ and $t \in [0, 1]$) is called a *controlled Markov process* associated with $u^n$. The generator of the controlled process is given by

$$(\bar{\mathcal{L}}^n f)(y,t) \doteq n \sum_{v \in \mathbb{Z}^N} u^n(y,v,t)[f(y+v/n) - f(y)],$$

where $f : \mathcal{S}^n \to \mathbb{R}$ is continuous and bounded. The dependence of $\bar{\mathcal{L}}^n$ on $u^n$ is omitted in our notation. The original queueing model corresponds to the particular choice $u^n(y,v,t) = r(y,v)$.

The controls that are considered are not required to be bounded, and in fact it will be necessary to consider unbounded controls. Existence of controlled processes is therefore not automatic. Let $\mathcal{Q}$ denote the set of all bounded functions $\varphi : \mathcal{S}^n \times [0,1] \to \mathbb{R}$ for which $\varphi(y,t)$ is continuously differentiable in $t$. We say that the control $u^n$ *has an associated controlled process* if there exists a Markov process $\{\Xi^n(t), t \in [0,1]\}$ on some probability space such that $\Xi^n(0) = y$ with probability 1 and such that

$$\varphi(\Xi^n(t),t) - \varphi(y,0) - \int_0^t \left( \frac{\partial \varphi(\Xi^n(s),s)}{\partial s} + \bar{\mathcal{L}}^n \varphi(\Xi^n(s),s) \right) ds \qquad (2)$$

is a martingale in $t$ for $t \in [0,1]$ for all $\phi \in \mathcal{Q}$, and if all processes satisfying these two conditions have the same distribution. In particular, every bounded control with a finite number of jump directions has an associated controlled process. We refer to a control $u^n$ as an *admissible* control if it satisfies (1) and has an associated controlled process.

For $\varphi \in D([0,1] : \mathcal{S}^n)$, $\beta \in \mathbb{R}^N$, and $\epsilon > 0$, let the *infinite exit cost* be defined as

$$g(\varphi;\beta,\epsilon) \doteq \begin{cases} \infty & \text{if there is } t \in [0,1] \text{ such that } \|\varphi(t) - t\beta\| \geq \epsilon, \\ 0 & \text{otherwise.} \end{cases}$$

For $a \in \mathbb{R}$, let

$$\ell(a) \doteq \begin{cases} a \log a - a + 1 & \text{if } a \geq 0, \\ \infty & \text{otherwise,} \end{cases}$$

where by convention $0 \log 0 \doteq 0$ and $0\ell(0/0) = 0$. Note that $\ell$ is nonnegative and that it has superlinear growth, i.e.,

$$\lim_{a \to \infty} \ell(a)/a = \infty.$$

Let $\bar{E}_y^n$ denote expectation conditioned on $\Xi^n(0) = y$.

**Theorem 1 (Dupuis and Ellis [3])** *Assume Condition 1. Then for each $n \in \mathbb{N}, y \in \mathcal{S}^n, \beta \in \mathbb{R}^N$, and $\epsilon > 0$, one has*

$$q^n(y;\beta,\epsilon) = \inf \bar{E}_y^n \left\{ \int_0^1 \sum_{v \in \mathbb{Z}^N} r(\Xi^n(t),v) \ell\left( \frac{u^n(\Xi^n(t),v,t)}{r(\Xi^n(t),v)} \right) dt + g(\Xi^n;\beta,\epsilon) \right\}, \qquad (3)$$

*where the infimum is taken over all admissible controls $u^n$ and associated controlled processes $\Xi^n$.*

We use the term *running cost* for the integral in (3), namely for

$$\int_0^1 \sum_{v \in \mathbb{Z}^N} r(\Xi^n(t),v) \ell\left( \frac{u^n(\Xi^n(t),v,t)}{r(\Xi^n(t),v)} \right) dt,$$

and the term *expected running cost* for its expectation under $\bar{P}_y^n$. By the *exit cost* and the *expected exit cost* we refer to $g(\Xi^n; \beta, \epsilon)$ and to its expectation under $\bar{P}_y^n$, respectively.

We define also a finite exit cost and attach to it a control problem similar to the one considered in Theorem 1, the only difference being in the exit cost. For $\varphi \in D([0,1]: \mathcal{S}^n)$, $\beta \in \mathbb{R}^N$, $\epsilon > 0$ and $M \geq 0$, let

$$g_M(\varphi; \beta, \epsilon) \doteq g(\varphi; \beta, \epsilon) \wedge M,$$

and define

$$q_M^n(y; \beta, \epsilon) \doteq \inf \bar{E}_y^n \left\{ \int_0^1 \sum_{v \in \mathbb{Z}^N} r(\Xi^n(t), v) \ell \left( \frac{u^n(\Xi^n(t), v, t)}{r(\Xi^n(t), v)} \right) dt + g_M(\Xi^n; \beta, \epsilon) \right\}, \qquad (4)$$

where the infimum is again taken over all admissible controls $u^n$ and associated controlled processes $\Xi^n$.

The following condition on the model is referred to in [3] as the Communication/Controllability Condition. Although a weaker condition is also stated in [3] under which the LDP is proved, the condition below suffices to cover the models that we consider. If this condition holds for the full model, then it will hold for all local models that are obtained from it as well (under Condition 1).

**Condition 2** *There exists a number $K_0$ such that for each pair of points $x$ and $y$ in $\mathcal{S}$ there exists $J \in \mathbb{N}$ satisfying $J \leq K_0 \|x - y\|$, and a sequence of points $\{x_0, x_1, \ldots, x_J\}$ in $\mathcal{S}$, for which $x_0 = x$, $x_J = y$, and $r(x_j, x_{j+1} - x_j) > 0$ for all $j = 0, 1, \ldots, J-1$.*

The next theorem asserts that certain large deviation limits exist for a local model with state space $\mathbb{Z}_+^{N,K}$ if one considers a constant velocity trajectory for which the velocity $\beta$ lies in $\mathcal{F}^{K,K}$. Thus if the state space imposes a non-negativity constraint $\langle x, e_i \rangle = 0$, then $\langle \beta, e_i \rangle = 0$. Note that such a trajectory will lie in the closure of all facets of the local model. As shown in [3], these turn out to be the only velocities needed in order to define the rate function at the process level. The first part of the theorem is taken from [3]. The second part concerns the related finite exit cost problem, and its proof is deferred to the Appendix.

**Theorem 2 (Dupuis and Ellis [3])** *Let $\mathcal{S}$ and $r$ be given and assume Conditions 1 and 2 are satisfied. Let $K$ be as in Condition 1, and let $\beta \in \mathcal{F}^{K,K}$ be given. Then there exists a number $L(\beta) \in [0, \infty)$ such that the following holds.*

1.

$$\lim_{\epsilon \to 0} \lim_{\delta \to 0} \liminf_{n \to \infty} \inf_{\{y \in \mathcal{S}^n : \|y\| \leq \delta\}} q^n(y; \beta, \epsilon)$$
$$= \lim_{\epsilon \to 0} \lim_{\delta \to 0} \limsup_{n \to \infty} \sup_{\{y \in \mathcal{S}^n : \|y\| \leq \delta\}} q^n(y; \beta, \epsilon)$$
$$= L(\beta),$$

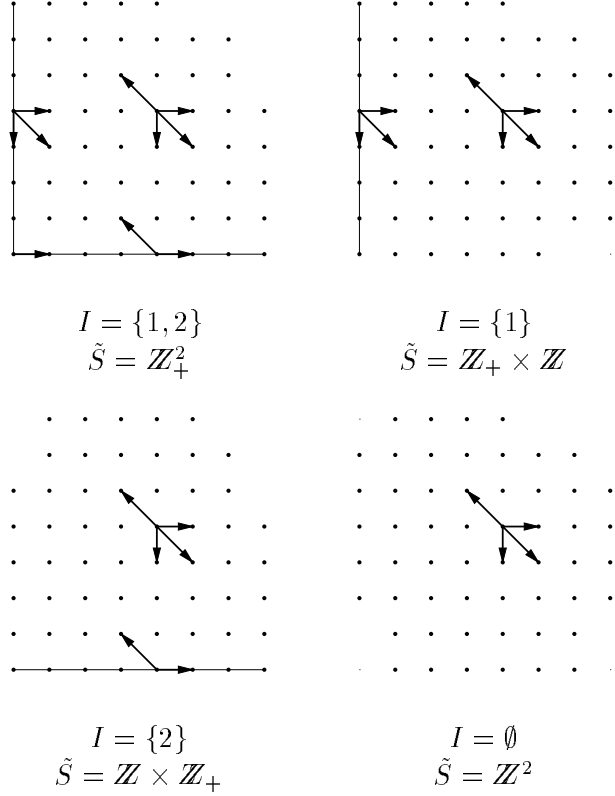

Figure 3: The state spaces $\tilde{\mathcal{S}}$ and intensity functions $\tilde{r}$ for different values of $I$ in a two dimensional example

2.
$$\lim_{M \to \infty} \lim_{\epsilon \to 0} \lim_{\delta \to 0} \liminf_{n \to \infty} \inf_{\{y \in \mathcal{S}^n : \|y\| \leq \delta\}} q_M^n(y; \beta, \epsilon)$$
$$= \lim_{M \to \infty} \lim_{\epsilon \to 0} \lim_{\delta \to 0} \limsup_{n \to \infty} \sup_{\{y \in \mathcal{S}^n : \|y\| \leq \delta\}} q_M^n(y; \beta, \epsilon)$$
$$= L(\beta).$$

**Remark:** Under the assumptions of Theorem 2, it is proved in [3] that the function $L(\cdot)$ is finite and convex on the linear subspace $\mathcal{F}^{K,K}$.

As discussed previously, the result of Theorem 2 on the local model serves as the basis for the full LDP. We next define the rate function for the full LDP in terms of the rate functions associated with local models. We therefore return to the full model with state space $\mathcal{S} = \mathbb{Z}_+^N$ and a given intensity function $r$, and assume for this full model that Conditions 1 and 2 are satisfied (with $K$ in those conditions equal to $\{1, 2, \ldots, N\}$). Consider a point $x$ in $\mathbb{R}_+^N$, and let $I = I(x) \doteq \{i : \langle x, e_i \rangle = 0\}$. Then the local model associated with such a point has state space $\tilde{\mathcal{S}} \doteq \mathbb{Z}_+^{N,I}$, and the only facet on which the local rate function need be defined is $\mathcal{F}^{I,I}$. The correct localized version of the jump intensities is defined for any $v \in \mathbb{Z}^N$ as follows. If $\tilde{x}$ is in one of the possible facets $\mathcal{F}^{I,M}$, $M \subset I$ of the local model, then we associate a point that is in the intersection of this facet with the original state space, but away from all non-negativity constraints save those indexed by $I$. In other words, $x \in \mathcal{F}^{\{1,\ldots,N\},M} \cap \mathcal{S}$. We then set $\tilde{r}(\tilde{x}, v) = r(x, v)$. Loosely speaking, in the manner suggested by Example 1 we throw away the parts of the generator associated parts of the domain where constraints other than those in $I$ are active, and in a corresponding manner extend the state space. A simple two dimensional example is illustrated in Figure 3 for the case of a Jackson network. We automatically get that the assumptions of Theorem 2 are satisfied by each possible $(\tilde{\mathcal{S}}, \tilde{r})$, and therefore conclude for each local model the existence of a function $\tilde{L}_I(\cdot)$ defined on $\mathcal{F}^{I,I}$. We then define $L(x, \beta)$ for any $x \in \mathbb{R}_+^N$ and $\beta \in \mathbb{R}^N$ by $L(x, \beta) = \tilde{L}_I(\beta)$ whenever $I = I(x)$ and $\beta \in \mathcal{F}^{I,I}$, and $L(x, \beta) = \infty$ otherwise.

For $S = \mathbb{R}^N$ or $S = \mathbb{R}_+^N$, let $\mathcal{T}([0,1] : S)$ denote the subset of $D([0,1] : S)$ of piecewise linear functions whose derivative has finitely many discontinuities. For $x \in \mathbb{R}_+^N$ and $\phi \in \mathcal{T}([0,1] : \mathbb{R}_+^N)$ satisfying $\phi(0) = x$, let

$$\tilde{J}_x(\phi) \doteq \int_0^1 L(\phi(t), \dot{\phi}(t)) dt.$$

Since $\phi \in \mathcal{T}([0,1] : \mathbb{R}_+^N)$, we have that $\dot{\phi}(t)$ is well defined, and moreover $\dot{\phi}(t) \in \mathcal{F}^{I(\phi(t)),I(\phi(t))}$ on the intervals where $\dot{\phi}$ is constant. Hence $L(\phi(t), \dot{\phi}(t)) < \infty$ for all but finitely many points $t$, and therefore $\tilde{J}_x(\phi)$ is well defined and finite. For all other $x \in \mathbb{R}^N$ and $\phi \in D([0,1], \mathbb{R}^N)$, we set $\tilde{J}_x(\phi) = \infty$.

Let $\sigma(\cdot, \cdot)$ denote the Skorokhod metric on $D([0,1] : \mathbb{R}^N)$, and for $\psi \in D([0,1] : \mathbb{R}^N)$ and $\alpha > 0$ let

$$B_\sigma(\psi, \alpha) \doteq \{\zeta \in D([0,1] : \mathbb{R}^N) : \sigma(\zeta, \psi) < \alpha\}.$$

Finally, the lower semicontinuous regularization of $\tilde{J}_x$ is denoted by $J_x$, namely for $x \in \mathbb{R}^N$ and $\psi \in D([0,1] : \mathbb{R}^N)$,

$$J_x(\psi) \doteq \lim_{\alpha \to 0} \inf_{y \in \mathbb{R}^N : \|y-x\| < \alpha} \inf_{\phi \in B_\sigma(\psi, \alpha)} \tilde{J}_y(\phi).$$

**Theorem 3 (Dupuis and Ellis [3])** *Assume Conditions 1 and 2. Let $x \in \mathbb{R}_+^N$. Then the following conclusions hold.*

1. *The function $J_x$ is nonnegative and lower semicontinuous on $D([0,1] : \mathbb{R}^N)$. Furthermore, for any compact set $C \subset \mathbb{R}^N$ and any $M \in [0, \infty)$, the set*

$$\bigcup_{x \in C} \{\psi \in D([0,1] : \mathbb{R}^N) : J_x(\psi) \leq M\}$$

   *is compact in $D([0,1] : \mathbb{R}^N)$.*

2. *For any open set $G$ in $D([0,1] : \mathbb{R}^N)$ and each $x \in \mathbb{R}^N$ we have the large deviation lower bound*

$$\lim_{\delta \to 0} \liminf_{n \to \infty} \inf_{y \in \mathcal{S}^n : \|y-x\| \leq \delta} \frac{1}{n} \log P_y^n(X^n \in G) \geq -\inf_{\psi \in G} J_x(\psi).$$

3. *For any closed set $F$ in $D([0,1] : \mathbb{R}^N)$ and each $x \in \mathbb{R}^N$ we have the large deviation upper bound*

$$\lim_{\delta \to 0} \limsup_{n \to \infty} \sup_{y \in \mathcal{S}^n : \|y-x\| \leq \delta} \frac{1}{n} \log P_y^n(X^n \in F) \leq -\inf_{\psi \in F} J_x(\psi).$$

With the preceding result available, all that remains is to identify the rate function $J_x$. Although it can easily be identified in certain special cases (e.g., for stable 2-dimensional models), this identification is in general a rather difficult problem. In the next two sections we introduce tools that will allow us to precisely characterize $L(x,\beta)$, and then apply them in Section 5 to two interesting classes of multidimensional models. This will identify the function $\tilde{J}_x$ on $\mathcal{T}([0,1]:\mathbb{R}^N_+)$. As we will see, a necessary condition for these methods to work is that the upper large deviation rate function obtained in [6] must equal $J_x$ on $\mathcal{T}([0,1]:\mathbb{R}^N_+)$. This fact and regularity properties of the various rate functions will be used in Section 4 to show that $J_x$ takes the expected form

$$J_x(\phi) \doteq \int_0^1 L(\phi(t), \dot{\phi}(t))dt$$

for all absolutely continuous functions $\phi$.

## 2.2 The Skorokhod Problem

In this subsection we give the precise definition of the Skorokhod Problem. The regularity properties of this problem play a key role in characterizing the limits of the controlled processes $\Xi^n$, which are in turn used to evaluate the limits of the quantities $q^n$ that appear in Theorem 2.

For a closed set $G \subset \mathbb{R}^N$, let

$$D_G([0,1]:\mathbb{R}^N) \doteq \{\psi \in D([0,1]:\mathbb{R}^N) : \psi(0) \in G\}.$$

For each point $x$ on the boundary of $G$ one is given a set $d(x)$ of unit vectors in $\mathbb{R}^N$. The Skorokhod Map assigns to every path $\psi \in D_G([0,1]:\mathbb{R}^N)$ a path $\phi$ that starts at $\phi(0) = \psi(0)$, but is constrained to $G$, in such a way that whenever it is in the interior of $G$ it is obtained by a translation of $\psi$, while on the boundary an additional "force" is used, allowed only in the directions defined by $d(x)$, so as to keep the path inside $G$. The definition of the Skorokhod Problem is stated below. For $\eta \in D_G([0,1]:\mathbb{R}^N)$ and $t \in [0,1]$ we let $|\eta|(t)$ denote the total variation of $\eta$ on $[0,t]$ with respect to the Euclidean norm on $\mathbb{R}^N$.

**Definition 1** *Let $\psi \in D_G([0,1]:\mathbb{R}^N)$ be given. Then $(\phi, \eta)$ solves the SP for $\psi$ with respect to $G$ and $d$ if $\phi(0) = \psi(0)$, and if for all $t \in [0,1]$ one has*

1. *$\phi(t) = \psi(t) + \eta(t)$,*
2. *$\phi(t) \in G$,*
3. *$|\eta|(t) < \infty$,*
4. *$|\eta|(t) = \int_{[0,t]} \mathbf{1}\{\phi(\theta) \in \partial G\} d|\eta|(\theta)$,*
5. *There exists a Borel measurable function $\gamma : [0,1] \to \mathbb{R}^N$ such that $d|\eta|$-almost everywhere one has $\gamma(t) \in d(\phi(t))$, and such that*

$$\eta(t) = \int_{[0,t]} \gamma(\theta) d|\eta|(\theta).$$

Note that $\eta$ changes only when $\phi$ is on the boundary, and only in the directions $d(\phi)$. The function $\Gamma$ defined by $\phi = \Gamma(\psi)$, on the domain where there is a unique solution to the SP, is called the Skorokhod Map. The functions $\psi$, $\phi$ and $\eta$ are referred to as the unconstrained, the constrained and the constraining paths, respectively.

In this paper we shall consider only sets $G$ of the following form:

$$G \doteq \bigcap_{i=1}^{q} \{x \in \mathbb{R}^N : \langle x, n_i \rangle \geq 0\}, \tag{5}$$

for some finite set of unit vectors $\{n_i, i = 1, 2, \ldots, q\}$. To any vector $n_i$ we attach a unit vector $d_i$ such that $\langle d_i, n_i \rangle > 0$. For $x \in \partial G$ let $I(x) \doteq \{i : \langle x, n_i \rangle = 0\}$. For any $x \in \partial G$ we define

$$d(x) \doteq \left\{ \gamma = \sum_{i \in I(x)} \alpha_i d_i : \alpha_i \geq 0, \|\gamma\| = 1 \right\}. \tag{6}$$

It is not assumed that $\{n_i, i = 1, 2, \ldots, q\}$ provides a minimal representation for $G$, and in fact this redundancy is necessary in some cases (including one of our examples) in order to allow for the proper full set of directions $d(x)$.

Both the domain $G$ and the sets of directions $d(x)$ are now entirely defined by $\{(n_i, d_i) : i = 1, 2, \ldots, q\}$, and therefore so is the SP.

**Definition 2** *A Skorokhod Map $\Gamma$ is called regular if it is defined on a subset $\mathcal{G}$ of $D_G([0,1] : \mathbb{R}^N)$ that includes all paths of bounded variation, and if it is Lipschitz continuous in that there is $K_1 < \infty$ such that for all $\psi_1$ and $\psi_2$ in $\mathcal{G}$*

$$\sup_{t \in [0,1]} \|\Gamma(\psi_1)(t) - \Gamma(\psi_2)(t)\| \leq K_1 \sup_{t \in [0,1]} \|\psi_1(t) - \psi_2(t)\|.$$

A SP is called regular if the corresponding SM is regular. If a SM is regular then there exists a unique Lipschitz continuous extension of $\Gamma$ to all of $D_G([0,1] : \mathbb{R}^N)$ [7], and with an abuse of notation this extension will also be denoted by $\Gamma$.

## 3 Representations for the Local Rate Function

As discussed in Section 2, local models will be used to calculate the rate function for a full queueing model. The local models simplify the problem by eliminating from consideration those parts of the process that are not involved in determining the rate function $L(x, \beta)$ at a point $x$. However, in all cases (save the local model associated with the interior of the full model) we are still obliged to deal with processes whose generator is discontinuous in the state variable.

As discussed previously, if one wishes to consider the probability that the rescaled process $X^n$ stays near a trajectory that evolves in an $(N-m)$–dimensional facet, then one must consider how this likelihood is affected by the form of the generator in each nearby region of constant statistical behavior (i.e., each facet whose closure contains the given $(N - m)$–dimensional facet). Typically, there will be $2^m$ regions of different statistical behavior that will be relevant for such a trajectory. A key quantity involved is the asymptotic fraction of time that the controlled process spends in each such region, and in general it is very hard to obtain this sort of information. However, when certain structural properties are present, one can identify alternative variables that provide all the information needed to identify the rate function, and yet which are easily obtained as a function the jump rates and jump vectors. These variables in fact arise from alternative representations of the local rate function, and the first results of this section will introduce these variables and indicate their connection with standard representations of the local rate function.

We therefore return to the "local model" setting of Theorem 1, and assume that a state space $\mathcal{S}$ and an intensity function $r$ are given, and that they satisfy Condition 1. The set $K$ appearing in Condition 1 will be arbitrary but fixed. According to Section 2.1, the rate function for the full model will be determined if we identify, for each such $K$, the function $L(\beta)$ that satisfies the conclusion of Theorem 2 for all $\beta \in \mathcal{F}^{K,K}$. As in Section 2.1, to simplify the notation we omit the dependence on $K$ from $\mathcal{S}, L(\beta), r(x, v)$, and so on.

We recall that if $x$ is in the facet $\mathcal{F}^{K,I}$ of the local model, then there are no constraints on $x_i$ if $i \notin K$, that $x_i \geq 0$ if $i \in K$, and that $x_i = 0$ if in addition $i \in I$. According to Condition 1, $r(x, v)$ is independent of $x$ for all $x$ belonging to a given facet $\mathcal{F}^{K,I}$, $I \subset K$. It is convenient to introduce the notation $r_{I,v} = r(x, v)$, where $x$ is any point in $\mathcal{F}^{K,I}$. Thus $r_{I,v}$ is the jump intensity from any location in the facet $\mathcal{F}^{K,I}$ in direction $v \in \mathbb{Z}^N$. Unless explicitly noted otherwise, all sums on $I$ in this section will be over all subsets of $K$, and all sums on $v$ will be over $V$.

For numbers that we denote by $\rho_I, u_{I,v}, \bar{r}_v$ and $c_v$, we shall consider the following sets of conditions:

$$\begin{cases} \rho_I \geq 0, \\ \sum_I \rho_I = 1, \\ u_{I,v} \geq 0, \\ r_{I,v} = 0 \implies u_{I,v} = 0, \end{cases} \tag{7}$$

$$\sum_{I,v} \rho_I u_{I,v} v = \beta. \tag{7a}$$

$$\begin{cases} \rho_I \geq 0, \\ \sum_I \rho_I = 1, \\ \bar{r}_v = \sum_I \rho_I r_{I,v}, \\ c_v \geq 0, \\ \sum_v \bar{r}_v c_v v = \beta. \end{cases} \tag{8}$$

Loosely speaking, these quantities may be interpreted as follows: $u_{I,v}$ will be the controlled jump rate in facet $\mathcal{F}^{K,I}$; $\rho_I$ represents the asymptotic fraction of time that the process which uses these controls spends in this facet, so that equation (7a) implies that the mean velocity of the controlled process is $\beta$; $\bar{r}_v$ is the average of the original jump rates in the direction $v$, with the averaging done according to the weights $\rho_I$; $c_v$ is a multiplier that is independent of $I$, and which determines the ration of $u_{I,v}/r_{I,v}$. Thus the second set of conditions represents a more structured system, since we essentially consider only controls with a certain type of independence from $I$.

Typically, there will be $N$ equality constraints on the $\rho_I$ from the velocity equation (7a), and one further equality constraint from the fact that the $\rho_I$'s are a probability vector. On the other hand, the number of unknown values of $\rho_I$ can be as large as $2^N$. This in part explains the simplicity of the stable two dimensional case, which one can in fact reduce to the case $N = 1$ [11]. In general, however, the $\rho_I$'s are not well defined.

In spite of this ambiguity, for some models one can use Jensen's inequality to restrict the class of controls under consideration, and in this restricted class replace the $\rho_I$'s by a weighted average, as suggested by the middle equality in (8). In such circumstances, the class of controlled models will be indexed by the collection of multipliers $\{c_v, v \in V\}$. Thus one can consider these multipliers as playing the role formerly played by the $u_{I,v}$. The significance of the parameters $\bar{r}_v$ will then follow from the equality of (9) and (10) below, since the independence of $c_v$ from $I$ means that the sum on $I$ in (7a) can be taken first, effectively replacing the quantities $\rho_I r_{I,v}$ (with $r_{I,v}$ known but $\rho_I$ unknown) by the quantities $\bar{r}_v$, and (7a) by the last equality in (8).

Now all this would be of little consequence if the $\bar{r}_v$ were as poorly defined as the $\rho_I$, but for the models we consider this is not true. In fact, it will turn out that the third and fifth equalities in (8) will provide the same number of equality relations as the cardinality of $V$.

We begin with a lower bound on the large deviations rate function for the local model, which corresponds to a large deviation upper bound.

**Lemma 1** *Under Conditions 1 and 2*

$$L(\beta) \geq \inf\left\{\sum_{I,v} \rho_I r_{I,v} \ell(u_{I,v}/r_{I,v}) : \text{Numbers } \rho_I \text{ and } u_{I,v} \text{ satisfying (7) and (7a)}\right\} \quad (9)$$

$$= \inf\left\{\sum_v \bar{r}_v \ell(c_v) : \text{Numbers } \rho_I, \bar{r}_v \text{ and } c_v \text{ satisfying (8)}\right\} \quad (10)$$

**Remark:** The bound (9), though with rather different notation, has already been established in [6].

**Proof:** We consider the following relaxed version of (7a):

$$\left\|y + \sum_{I,v} \rho_I u_{I,v} v - \beta\right\| < \epsilon. \quad (7b)$$

It will be shown below that for any $n \in \mathbb{N}$, $y \in \mathcal{S}^n$, $\beta \in \mathbb{R}^N$, $\delta > 0$ and $\epsilon > 0$, such that $\|y\| \leq \delta < \epsilon$,

$$q^n(y; \beta, \epsilon) \geq \inf\left\{\sum_{I,v} \rho_I r_{I,v} \ell(u_{I,v}/r_{I,v}) : \text{Numbers } \rho_I \text{ and } u_{I,v} \text{ satisfying (7) and (7b)}\right\}. \quad (11)$$

Given this inequality, part 1 of Theorem 2 and the lower semicontinuity of $\ell$ imply (9). Moreover, if the numbers $\rho_I$, $\bar{r}_v$ and $c_v$ satisfy (8) and if $u_{I,v} = c_v r_{I,v}$, then $\rho_I$ and $u_{I,v}$ satisfy (7) and (7a). Therefore the right hand side of (9) is less than or equal to the expression in (10). On the other hand, for any $\rho_I$ and $u_{I,v}$ satisfying (7) and (7a), the numbers $\rho_I$, $\bar{r}_v \doteq \sum_I \rho_I r_{I,v}$ and $c_v \doteq \sum_I \rho_I u_{I,v} / \sum_I \rho_I r_{I,v}$ satisfy (8). Jensen's inequality applied to the convex function $\ell$ implies that

$$\sum_{I,v} \rho_I r_{I,v} \ell(u_{I,v}/r_{I,v}) \geq \sum_v \bar{r}_v \ell(c_v),$$

and thus (10) follows. We therefore turn to the proof of (11), relying on the representation (3).

Obviously, the exit cost in (3) may be omitted if we infimize over all controls which make it vanish, namely

$$q^n(y; \beta, \epsilon) = \inf \bar{E}_y^n \left\{ \int_0^1 \sum_{v \in \mathbb{Z}^N} r(\Xi^n(t), v) \ell\left( \frac{u^n(\Xi^n(t), v, t)}{r(\Xi^n(t), v)} \right) dt \right\}, \quad (12)$$

where the infimum is taken over the set of admissible controls $u^n$ for which the associated controlled processes $\Xi^n$ satisfy

$$\bar{P}_y^n(g(\Xi^n; \beta, \epsilon) = 0) = 1. \quad (13)$$

We use the following notation:

$$\rho_y^n(I) \doteq \bar{E}_y^n \int_0^1 \mathbf{1}\{\Xi^n(t) \in \mathcal{F}^{K,I}\} dt, \text{ and } U_y^n(I, v) \doteq \frac{1}{\rho_y^n(I)} \bar{E}_y^n \int_0^1 u(\Xi^n(t), v, t) \mathbf{1}\{\Xi^n(t) \in \mathcal{F}^{K,I}\} dt \quad (14)$$

if $\rho_y^n(I) > 0$, while $U_y^n(I, v) \doteq 0$ if $\rho_y^n(I) = 0$. Suppose that we rewrite the expected running cost in (12) as

$$\bar{E}_y^n \left\{ \sum_{I,v} \int_0^1 \mathbf{1}\{\Xi^n(t) \in \mathcal{F}^{K,I}\} r_{I,v} \ell\left( \frac{u^n(\Xi^n(t), v, t)}{r_{I,v}} \right) dt \right\}.$$

It then follows from Jensen's inequality [applied to the convex function $\ell(\cdot)$] that (12) is greater than or equal to

$$\sum_{I,v} \rho_y^n(I) r_{I,v} \ell(U_y^n(I, v)/r_{I,v}).$$

Since $u^n$ is an admissible control, we know that $r_{I,v} = 0$ implies $U_y^n(I, v) = 0$, and also that for any $M > 0$ and $j \in \{1, 2, \ldots, N\}$ the process defined in (2) is a martingale if $\varphi(x, t) = \langle x, e_j \rangle \wedge M$. Taking $M \to \infty$, one obtains by the monotone convergence theorem that

$$\bar{E}_y^n \Xi^n(1) = y + \bar{E}_y^n \int_0^1 \sum_v u^n(\Xi^n(t), v, t) v = y + \sum_{I,v} \rho_y^n(I) U_y^n(I, v) v.$$

The constraint (13) implies that $\rho_I \doteq \rho_y^n(I)$ and $u_{I,v} \doteq U_y^n(I, v)$ satisfy (7) and (7b). This proves (11), and completes the proof of the lemma.
□

We next derive an upper bound on the rate function $L(\beta)$. The proof of the lower bound just given suggests that it may be sufficient to consider controls that are constant in the time variable. If such a control is used, and if in addition it is of the special form associated with the set of constraints (8) (i.e., $u^n(x,v,t) = c_v r_{I,v}$ for $x \in \mathcal{F}^{K,I}$, $I \subset K$, $v \in \mathbb{Z}^N$, and $t \in [0,1]$), then the definition (4) of $q_M^n$ implies that

$$q_M^n(y;\beta,\epsilon) \leq \sum_{I,v} \rho_y^n(I) r_{I,v} \ell(c_v) + M \bar{P}_y^n(\sup_t \|\Xi^n(t) - \beta t\| \geq \epsilon) \qquad (15)$$

for any such admissible control $u^n$.

We introduce two conditions that will be needed in the treatment of the upper bound. Both parts refer to (8). The first states that the collection $\{\bar{r}_v\}$, when viewed as a function of the $\{c_v\}$ and $\beta \in \mathcal{F}^{K,K}$, is unique (if it exists) and continuous. When this property holds it identifies the $\bar{r}_v$ as a quantity that is easier to work with (when compared to the $\rho_I$'s). The second part of the condition is less significant, and simply requires that for any solution $c_v \geq 0$, $v \in V$, there exist a nearby solution with $c_v > 0$, $v \in V$.

**Condition 3** *We consider $\mathcal{S}$ and $r$ satisfying Condition 1 and let $K$ be such that $\mathcal{S} = \mathbb{Z}_+^{N,K}$.*

1. *Let $c_v > 0, v \in V$ be given. Then for all $\epsilon > 0$ there exists $\delta > 0$ such that if $(\{\rho_I\}_{I \subset K}, \{\bar{r}_v\}_{v \in V}, \{c_v\}_{v \in V}, \beta)$ and $(\{\rho'_I\}_{I \subset K}, \{\bar{r}'_v\}_{v \in V}, \{c_v\}_{v \in V}, \beta')$ satisfy (8), where $\beta, \beta' \in \mathcal{F}^{K,K}$ and $\|\beta - \beta'\| < \delta$, then $|\bar{r}_v - \bar{r}'_v| < \epsilon$, $v \in V$.*

2. *Let the quadruple $(\{\rho_I\}_{I \subset K}, \{\bar{r}_v\}_{v \in V}, \{c_v\}_{v \in V}, \beta)$, where $c_v \geq 0$, $v \in V$ and $\beta \in \mathcal{F}^{K,K}$, satisfy (8). Then there exists an $M < \infty$ such that for all $\epsilon > 0$ there exists $(\{\rho'_I\}_{I \subset K}, \{\bar{r}'_v\}_{v \in V}, \{c'_v\}_{v \in V}, \beta')$ satisfying (8), such that*

    - $c'_v > 0$, $v \in V$, $\beta' \in \mathcal{F}^{K,K}$, $|\bar{r}_v - \bar{r}'_v| < \epsilon$, $v \in V$, and $\|\beta - \beta'\| < \epsilon$,
    - *for $v \in \{w \in V : \bar{r}_w \neq 0\}$ $|c_v - c'_v| < \epsilon$, and for $v \in \{w \in V : \bar{r}_w = 0\}$ $c'_v \leq M$.*

Note that part 1 of Condition 3 implies uniqueness of the numbers $\{\bar{r}_v\}$ given $\beta$ and $c_v > 0$. However, existence of $\bar{r}_v$ is not assumed. Also, we reiterate that it may happen (and indeed will happen in our applications) that there is uniqueness of the $\{\bar{r}_v\}$, without uniqueness of the corresponding $\{\rho_I\}$ appearing in (8).

Our proof of the lower bound is based on weak convergence of the trajectories of the controlled processes to the solution of a related SP. The following condition asserts that this SP is well behaved. The first part simply asserts that the domain of the SP and that of the large deviation problem are compatible, while the second states the nonnegativity relation between normals and directions of constraint that is needed in the SP. To identify the LLN limit of the controlled queueing network as the solution of a SP, we will represent the system as an average "drift" perturbed by correction terms. The fourth part ensures that these correction terms point in a direction that is consistent with the given SP. The key part is 3, which will identify the solution to the SP as the unique LLN limit of the controlled queueing system.

**Condition 4** *Consider $\mathcal{S}$ and $r$ satisfying Condition 1 and let $K$ be such that $\mathcal{S} = \mathbb{Z}_+^{N,K}$. Then for every set of numbers $c_v > 0$, $v \in V$, for which (8) holds with some $\beta \in \mathcal{F}^{K,K}$, $\{\rho_I\}$ and $\{\bar{r}_v\}$, there exist a number $q \in \mathbb{N}$, and unit vectors $\{n_i, i = 1, 2, \ldots, q\}$ and $\{d_i, i = 1, 2, \ldots, q\}$ (that may depend on $\{c_v\}$), such that the following hold.*

1. $\mathbb{R}_+^{N,K} = \bigcap_{i=1}^q \{x \in \mathbb{R}^N : \langle x, n_i \rangle \geq 0\}$.

2. $\langle d_i, n_i \rangle > 0$ for $i = 1, 2, \ldots, q$.

3. The SM associated with $\{(n_i, d_i), i = 1, 2, \ldots, q\}$ is regular in the sense of Section 2.2.

4. For any $I \subset K$, $I \neq \emptyset$, if $\sum_v c_v(r_{I,v} - r_{\emptyset,v})v \neq 0$ then

$$\frac{\sum_v c_v(r_{I,v} - r_{\emptyset,v})v}{\|\sum_v c_v(r_{I,v} - r_{\emptyset,v})v\|} \in d(x), \quad x \in \mathcal{F}^{K,I},$$

where $d(\cdot)$ is defined by (6) with the $d_i$ as given above.

Note that we do not assume the representation in part 1 is a minimal representation for $\mathbb{R}_+^{N,K}$.

We have the following.

**Proposition 1** *Consider $\mathcal{S}$ and $r$ satisfying Conditions 1 and 4, and let $K$ be such that $\mathcal{S} = \mathbb{Z}_+^{N,K}$. Let $c_v > 0$, $v \in V$ and $\beta \in \mathcal{F}^{K,K}$ be given, and assume (8) is satisfied with some $\{\rho_I\}$ and $\{\bar{r}_v\}$. For $n \in \mathbb{N}$ let $u^n$ be an admissible control defined by $u^n(x, v, t) = c_v r_{I,v}$ for $x \in \mathcal{F}^{K,I}$, $I \subset K$, $v \in V$ and $t \in [0, 1]$. Then for all $\epsilon > 0$,*

$$\lim_{n \to \infty} \bar{P}_0^n \left( \sup_t \|\Xi^n(t) - \beta t\| \geq \epsilon \right) = 0. \tag{16}$$

Note that according to our notation $\bar{P}_0^n(\Xi^n(0) = 0) = 1$.

The proof of this proposition is given in Section 4. Let us show that an upper bound on the rate function follows, under all of Conditions 1, 2, 3 and 4.

Let $\beta \in \mathcal{F}^{K,K}$ be fixed. Let $c_v > 0$, $v \in V$ be such that there exist $\{\rho_I\}$ and $\{\bar{r}_v\}$ for which (8) holds with the above $\beta$. Note that the $\{\bar{r}_v\}$ are unique, by part 1 of Condition 3. Furthermore, recall the notation in (14), and observe that in terms of this notation

$$\bar{E}_0^n \Xi^n(1) = \sum_{I,v} \rho_0^n(I) c_v r_{I,v} v.$$

It is easy to verify the bound $\sup_n \bar{E}_0^n \|\Xi^n(1)\|^2 < \infty$. Therefore by (16)

$$\limsup_{n \to \infty} \left\| \sum_{I,v} \rho_0^n(I) c_v r_{I,v} v - \beta \right\| \leq \epsilon.$$

Part 1 of Condition 3 then implies that for all $v \in V$

$$\limsup_{n \to \infty} \left| \sum_I \rho_0^n(I) r_{I,v} - \bar{r}_v \right| \leq \alpha_\epsilon,$$

where $\alpha_\epsilon \to 0$ as $\epsilon \to 0$. By part 2 of Theorem 2 and (15) we obtain

$$L(\beta) \leq \liminf_{M \to \infty} \liminf_{\epsilon \to 0} \liminf_{n \to \infty} q_M^n(0; \beta, \epsilon) \leq \sum_v \bar{r}_v \ell(c_v).$$

One may now take infimum over all $c_v > 0$, $v \in V$, $\{\rho_I\}$, $I \subset K$ and $\{\bar{r}_v\}$, $v \in V$ for which (8) holds. The bound that one obtains still differs from the expression in (10), where all the numbers $\{c_v\}$ are allowed to be zero. We now show that the two infima are the same. Suppose that we are given quantities that satisfy the conditions (8), with $c_v \geq 0$, $v \in V$. For $\epsilon > 0$ let the primed versions be associated as in part 2 of Condition 3. Hence $c'_v > 0$, $v \in V$, and using the inequality just proved and the bounds $|\bar{r}_v - \bar{r}'_v| < \epsilon$, $c'_v \leq M$, $v \in \{w \in V : \bar{r}_w = 0\}$, we obtain

$$L(\beta) \leq L(\beta) - L(\beta') + \sum_v \bar{r}'_v \ell(c'_v) \leq L(\beta) - L(\beta') + \sum_{v \in V : \bar{r}_v \neq 0} \bar{r}'_v \ell(c'_v) + \#\{v \in V : \bar{r}_v = 0\}\epsilon[1 \vee \ell(M)].$$

We also have that $\beta' \in \mathcal{F}^{K,K}$, $|c_v - c'_v| < \epsilon$ for $v \in \{w \in V : \bar{r}_w \neq 0\}$, $|\bar{r}_v - \bar{r}'_v| < \epsilon$ and $\|\beta - \beta'\| < \epsilon$. Recall that $\ell(a)$ is continuous for $a \in \mathbb{R}_+$, and $L(\beta)$ is lower semi-continuous for $\beta \in \mathcal{F}^{K,K}$. Sending $\epsilon \to 0$ we obtain the desired inequality with the unprimed quantities. When combined with Lemma 1 for the lower bound we obtain the following, which is the main result of this paper.

**Theorem 4** *Let $\mathcal{S}$ and $r$ be fixed, assume Conditions 1, 2, 3 and 4, and let $K$ be such that $\mathcal{S} = \mathbb{Z}_+^{N,K}$. Then the rate function for the local model is given by*

$$L(\beta) = \inf\left\{\sum_v \bar{r}_v \ell(c_v) : \text{Numbers } \rho_I, \bar{r}_v \text{ and } c_v \text{ satisfying (8)}\right\}. \tag{17}$$

As discussed at the end of Section 2.1, the rate function for the full model (on path space) is obtained as follows. We assume Conditions 1, 2, 3 and 4 for the full model. For each point $x \in \mathbb{R}_+^N$ we let $I = I(x) \doteq \{i : \langle x, e_i \rangle = 0\}$, and then associate to this point the proper local model. As we will see, Conditions 1, 2, 3 and 4 then hold for the local model, and Theorem 4 identifies the function $L(x, \beta)$ as $L_I(\beta)$ for $\beta \in \mathcal{F}^{I,I}$. We let $L(x, \beta) = \infty$ in all other cases. By Theorem 3 the large deviation principle for the sequence $\{X^n, n \in \mathbb{N}\}$ holds with the rate function $J_x(\phi)$, which is the lower semicontinuous regularization of the function $\tilde{J}_x$ defined on $\mathcal{T}([0,1] : \mathbb{R}_+^N)$ by

$$\tilde{J}_x(\phi) \doteq \int_0^1 L(\phi(t), \dot{\phi}(t))dt$$

if $\phi(0) = x$, and $\infty$ otherwise.

The following theorem provides an explicit formula for $J_x$.

**Theorem 5** *Assume Conditions 1, 2, 3 and 4, and define $L(x, \beta)$ and $J_x(\phi)$ as in the last paragraph. Then the following representation holds:*

$$J_x(\phi) = \int_0^1 L(\phi(t), \dot{\phi}(t))dt$$

*if $\phi(0) = x$ and if $\phi \in D([0,1] : \mathbb{R}_+^N)$ is absolutely continuous, while $J_x(\phi) = \infty$ in all other cases.*

**Proof:** The function $\tilde{J}_x(\phi)$ coincides on $\mathcal{T}([0,1] : \mathbb{R}_+^N)$ with the upper large deviation rate function obtained in [6]. In the notation of this paper, the upper rate function of [6] is just

$$\hat{J}_x(\phi) \doteq \int_0^1 L(\phi(t), \dot{\phi}(t))dt$$

if $\phi(0) = x$ and if $\phi \in D([0,1] : \mathbb{R}_+^N)$ is absolutely continuous, and $\hat{J}_x(\phi) = \infty$ otherwise. Since it is proved in [6] that
$$\left\{ \phi \in D([0,1] : \mathbb{R}_+^N) : \hat{J}_{\phi(0)}(\phi) \leq M, \phi(0) \in C \right\}$$
is compact for any compact set $C \subset \mathbb{R}^N$ and $M < \infty$, $\hat{J}_x$ is lower semicontinuous. Since $J_x$ is the lower semicontinuous regularization of a function that agrees with $\hat{J}_x$ on a dense subset of its domain of definition, $\hat{J}_x(\phi) \leq J_x(\phi)$ for all $D([0,1] : \mathbb{R}_+^N)$. To prove the reverse inequality, we must show that given $\phi \in D([0,1] : \mathbb{R}_+^N)$ and $\epsilon > 0$ there is $\phi_\epsilon \in \mathcal{T}([0,1] : \mathbb{R}_+^N)$ such that $\hat{J}_x(\phi_\epsilon) \leq J_x(\phi) + \epsilon$ and $\sup_{t \in [0,1]} \|\phi(t) - \phi_\epsilon(t)\| \leq \epsilon$. However, the construction of such a function can be carried out using an argument similar to the one used to prove Lemma 7.5.4 of [4], and is omitted.
□

In Section 5 we will verify all the assumptions of Theorem 4 for some interesting models. The following monotonicity result of the regularity of the SM shows that if part 4 of Condition 4 holds for the full model, then it holds for all corresponding local models as well (for a different version see [9]). The proof is given in the Appendix. It is easy to check that all remaining parts of this condition also hold once they are verified for the full model.

**Lemma 2** *Consider a SM on $G \subset \mathbb{R}^N$ associated with $\{(n_i, d_i), i = 1, \ldots, q\}$, where $\langle n_i, d_i \rangle > 0$ and $G$ is as in (5), and assume that the SM is regular. Let $x \in \partial G$ be fixed, and let $I(x) = \{i : \langle n_i, x \rangle = 0\}$. Then the SM associated with $\{(n_i, d_i), i \in I(x)\}$ is also regular.*

## 4 Weak convergence considerations

This section contains the proof of Proposition 1, which is based on showing that $\{\Xi^n\}$ converges in distribution to the solution of the SP defined in Condition 4. We use several ideas from [7], and in particular, the proof of Lemma 5 closely follows that of [7, Theorem 3.2].

We assume that a state space $\mathcal{S}$ and an intensity function $r$ satisfying Condition 1 are given. Moreover, letting $K$ be as in Condition 1, we are given also $\beta \in \mathcal{F}^{K,K}$. On $\mathcal{S}^n$, $n \in \mathbb{N}$ we consider controls of the form $u^n(x,v) = u_{I,v}$, $x \in I$, $I \subset K$ and $v \in V$, and jump Markov processes $\Xi^n(t)$, $t \in [0,1]$ starting at zero with probability 1, and with jump intensity $u^n(x,v)$ from $x$ to $x+v$. As in Section 3, the sum on $I$ will always be over subsets of $K$, and the sum on $v$ will be over $V$. We assume that for some $c_v > 0$, $v \in V$, the controlled jump rates are given by $u_{I,v} = c_v r_{I,v}$. We recall that $\emptyset$ indexes the interior facet relative to the local model, i.e., $\mathcal{F}^{K,\emptyset}$. Further, we are given a SP $\{(n_i, d_i), i = 1, \ldots, q\}$ that satisfies Condition 4.

The proof will make use of a jump Markov process $(\bar{X}^n(t), Y^n(t))$ on $\mathcal{S}^n \times (n^{-1} \mathbb{Z}^N)$. The process will start at $(0,0)$ with probability 1. The $Y^n$ component will be a homogeneous jump Markov process on $n^{-1} \mathbb{Z}^N$, whose generator will always be the same as that of $\Xi^n$ for points in the interior $\mathcal{F}^{K,\emptyset}$, while $\bar{X}^n$ will be equal in law to $\Xi^n$. Moreover, $\bar{X}^n$ and $Y^n$ will have identical increments when $\bar{X}^n(t) \in \mathcal{F}^{K,\emptyset}$, and independent increments when $\bar{X}^n(t) \notin \mathcal{F}^{K,\emptyset}$. For such a pair process the jump intensities are as follows. Let $\mathbf{x}' = \mathbf{x} + \mathbf{v}$. Then the jump intensity $u^n(\mathbf{x}, \mathbf{v})$ from $\mathbf{x} = (x, y)$ to $\mathbf{x}'$ is given by:

$$u^n(\mathbf{x}, \mathbf{v}) = \begin{cases} nu_{\emptyset,v} & \text{if either } x \in \mathcal{F}^{K,\emptyset} \text{ and } \mathbf{x}' = (x + n^{-1}v, y + n^{-1}v) \\ & \text{or } x \notin \mathcal{F}^{K,\emptyset} \text{ and } \mathbf{x}' = (x, y + n^{-1}v) \\ nu_{I,v} & \text{if } x \in \mathcal{F}^{K,I} \neq \mathcal{F}^{K,\emptyset} \text{ and } \mathbf{x}' = (x + n^{-1}v, y), \\ 0 & \text{otherwise.} \end{cases}$$

We define $Z^n(t) = \bar{X}^n(t) - Y^n(t)$, $t \in [0,1]$ and let $Z^n = \bar{Z}^n + \hat{Z}^n$ be the Doob-Meyer decomposition of $Z^n$, where $\bar{Z}^n$ is a process that is predictable on the filtration generated by $(\bar{X}^n, Y^n)$, and $\hat{Z}^n$ is a martingale on the same filtration. If one considers the process $Y^n$ as an unconstrained version of the controlled queueing system, and $\bar{X}^n$ as the "correct" constrained version (corresponding to $\psi$ and $\phi$ respectively in the formulation of the SP in Section 2), then $Z^n$ (and more precisely $\bar{Z}^n$) will play the role of the constraining term $\eta$ in the SP. This correspondence will turn out to be exact in the limit $n \to \infty$.

To prove this fact, we define

$$\gamma^n(t) \doteq \begin{cases} \frac{\sum_v (u(\bar{X}^n(t), v) - u_{\emptyset,v}) v}{\|\sum_v (u(\bar{X}^n(t), v) - u_{\emptyset,v}) v\|} & \text{if } \sum_v (u(\bar{X}^n(t), v) - u_{\emptyset,v}) v \neq 0, \\ 0 & \text{otherwise} \end{cases}$$

for $t \in [0,1]$, and let $R^n(t) \doteq |\bar{Z}^n|(t)$. Denoting the closed unit sphere in $\mathbb{R}^N$ by $B(0,1)$, we let

$$\Sigma' \doteq \mathbb{R}_+^{N,K} \times B(0,1),$$

and define measures $\mu^n$ on $\Sigma \doteq [0,1] \times \Sigma'$ by

$$\mu^n([0,t] \times A) = \int_{[0,t]} \mathbf{1}\{(\bar{X}^n(\theta), \gamma^n(\theta)) \in A\} dR^n(\theta).$$

Proposition 1 will turn out to be a consequence of the following three lemmas. For the notation used in the statements we refer to Condition 4.

**Lemma 3** *Let*

$$\beta_0 = \sum_v u_{\emptyset,v} v,$$

*and for $t \in [0,1]$ let*

$$\psi(t) = \beta_0 t, \quad \phi(t) = \beta t, \quad \text{and} \quad \eta(t) = (\beta - \beta_0) t.$$

*Then $(\phi, \eta)$ solve the SP for $\psi$ with respect to $\mathbb{R}_+^{N,K}$ and $\{(n_i, d_i), i = 1, ..., q\}$.*

**Lemma 4** *The family $\{(\bar{X}^n, Y^n, Z^n, \bar{Z}^n, R^n, \mu^n), n \in \mathbb{N}\}$ is tight.*

**Lemma 5** *Consider any subsequence of*

$$(\bar{X}^n, Y^n, Z^n, \bar{Z}^n, R^n, \mu^n),$$

*and let*

$$(\bar{X}, Y, Z, \bar{Z}, R, \mu)$$

*denote the limit of a weakly convergent subsubsequence. Then w.p.1 $(\bar{X}, Z)$ solves the SP for $Y$ with respect to $\mathbb{R}_+^{N,K}$ and $\{(n_i, d_i), i = 1, ..., q\}$.*

Before giving the proofs of these lemmas, we show that Proposition 1 follows. It follows from the standard functional law of large numbers that the weak limit of $Y^n$ is exactly $\psi$ of Lemma 3. Note that regularity of the SM implies uniqueness. Thus by Lemmas 4 and 5 and the usual argument by contradiction, the weak limit of $\bar{X}^n$ is just $\phi(t) = \beta t$, $t \in [0, 1]$, w.p.1. This implies (16), and Proposition 1 follows.
□

**Proof of Lemma 3:** Properties 1 and 3 of Definition 1 are obvious. Since Proposition 1 assumes $\beta \in \mathcal{F}^{K,K}$, it follows that $\phi(t) \in \mathcal{F}^{K,K}$ for all $t \in [0, 1]$. Since $\mathcal{F}^{K,K}$ is a subset of $\partial \mathbb{R}_+^{N,K}$, properties 2 and 4 follow. Next, if $\beta = \beta_0$ then $|\eta|(t) = 0$ for $t \in [0, 1]$, and property 5 obviously holds. If $\beta \neq \beta_0$, then let $\tilde{\gamma}(t) = (\beta - \beta_0)/\|\beta - \beta_0\|$ for $t \in [0, 1]$. It immediately follows that $\eta(t) = \int_0^t \tilde{\gamma}(\theta) d|\eta|(\theta)$. Since by (8)

$$\beta - \beta_0 = \sum_{I,v} c_v \rho_I (r_{I,v} - r_{\emptyset, v}) v,$$

it follows that $\tilde{\gamma}(t) \in d(x)$ if $x \in \mathcal{F}^{K,K}$ by part 4 of Condition 4. Since $\phi(t) \in \mathcal{F}^{K,K}$ for all $t \in [0, 1]$, property 5 of Definition 1 follows as well.
□

**Proof of Lemma 4:** That the family $\{(\bar{X}^n, Y^n), n \in \mathbb{N}\}$ is tight in the Skorokhod topology follows immediately from Aldous-Kurtz Theorem (see, e.g., [12]). Hence $\{Z^n\}$ is also tight. For the predictable part of $Z^n$ we have the following expression:

$$\bar{Z}^n(t) = \int_0^t \sum_v (u(\bar{X}^n(\theta), v) - u_{\emptyset, v}) v d\theta. \tag{18}$$

In particular, the trajectories of $\bar{Z}^n$ are Lipschitz with a common coefficient, and hence $\{\bar{Z}^n\}$ is tight. Similarly, $\{R^n\}$ is tight, and in fact there is $B < \infty$ such that $R^n(1) \leq B$ for all $n \in \mathbb{Z}$ w.p.1.

Thus all that remains is to show the tightness of $\{\mu^n\}$ as random variables in a space of measures with the weak topology. Let $C_\epsilon \subset D([0, 1] : \mathbb{R}_+^{N,K})$ be compact and such that

$$\bar{P}^n \{\bar{X}^n \in C_\epsilon\} \geq 1 - \epsilon$$

for all $\epsilon \in (0, 1)$ and all $n \in \mathbb{N}$. Then there is $M_\epsilon < \infty$ such that for any $\bar{X}^n \in C_\epsilon$

$$\|\bar{X}^n(t)\| \leq M_\epsilon < \infty \tag{19}$$

for all $t \in [0,1]$. We recall the bound $R^n(1) \leq B < \infty$, which implies $\mu^n(\Sigma) = R^n(1) \leq B$. If $\mu^n(\Sigma) \neq 0$ then $\mu^n/\mu^n(\Sigma)$ is a probability measure on $\Sigma$. Using (19), for all $\epsilon \in (0,1)$ there exists a compact set $C'_\epsilon \subset \Sigma$ such that for $n \in \mathbb{N}$

$$\bar{P}^n \left\{ \mu^n(\Sigma) \neq 0 \text{ and } \frac{\mu^n((C'_\epsilon)^c)}{\mu^n(\Sigma)} \geq \epsilon \right\} \leq \epsilon.$$

Hence

$$\bar{P}^n \left\{ \mu^n(\Sigma) \neq 0 \text{ and } \frac{\mu^n}{\mu^n(\Sigma)} \in C''_\delta \right\} \geq 1 - 2\delta, \tag{20}$$

where

$$C''_\delta = \bigcap_{m \in \mathbb{N}: 2^{-m} < \delta} \left\{ \nu \in \mathcal{P}(\Sigma) : \nu((C'_{2^{-m}})^c) \leq 2^{-m} \right\}.$$

For each $\delta \in (0,1)$ the set of probability measures $C''_\delta$ is tight, and therefore by Prohorov's Theorem it is also relatively compact. To complete the proof of tightness of the random measures $\mu^n$ we must also show that for all $\epsilon > 0$ there exists $M$ such that $P(\mu^n(\Sigma) > M) < \epsilon$. Since this is implied by the w.p.1 bound $\mu^n(\Sigma) \leq B$, the proof of the lemma is complete.

□

**Proof of Lemma 5:** According to the Skorokhod Representation Theorem, there exists a probability space on which there are random variables having distributions identical to $(X^n, Y^n, Z^n, \bar{Z}^n, R^n, \mu^n)$ and $(X, Y, Z, \bar{Z}, R, \mu)$, and for which the convergence is in the a.s. sense. In the proof of the lemma we will make use of this alternative space, but without changing the notation (see Shiryaev [16]).

In the proof to follow there will be numerous properties of and relations between the limit and prelimit random variables that hold only in a w.p.1 sense. To simplify the discussion, the w.p.1 qualifier will be omitted. Statements that hold only in an a.e. sense for the time variable will be explicitly identified, in which case the qualifier also holds w.p.1.

Suppose that the quadratic variation up to time $t$ of the martingale $\hat{Z}^n(t)$ is denoted by $\langle \hat{Z}^n \rangle(t)$. Then by the Burkholder-Davis-Gundy inequality, there exists a constant $c < \infty$ such that

$$E^n \sup_{t \in [0,1]} \|\hat{Z}^n(t)\|^2 \leq c E \langle \hat{Z}^n \rangle(1).$$

Since the jump rates and jump vectors are all uniformly bounded, there is $c_1 < \infty$ such that that $E\langle \hat{Z}^n \rangle(1) \leq c_1/n$. Therefore $\hat{Z}^n$ converges weakly to zero in the sup norm, and consequently $\bar{Z} = Z$.

We next observe that $\bar{X}^n(t) = Y^n(t) + Z^n(t)$ and also that $\bar{X}^n(t) \in \mathbb{R}_+^{N,K}$ for all $n \in \mathbb{N}$ and $t \in [0,1]$. The almost sure convergence implies the analogous equality and inclusion for the limit, and therefore properties 1 and 2 of Definition 1 follow. Recall from the proof of Lemma 4 that there is $B < \infty$ such that $R^n(1) \leq B$, which implies $R(1) \leq B$. Since the total variation of an element of $D([0,1]:S)$ (for any Polish space $S$) is a lower semi-continuous functional, we also have $|\bar{Z}|(1) < \infty$, and property 3 follows.

To prove properties 4 and 5 we first define the sets

$$\Sigma_1 = \{(t, x, \gamma) \in \Sigma : x \in \mathcal{F}^{K,\emptyset}\},$$

$$\Sigma_2 = \{(t, x, \gamma) \in \Sigma : \gamma \notin d(x)\}$$

and

$$\Sigma_3^\Delta = \{(t, x, \gamma) \in \Sigma : |x - X(t)| > \Delta\}.$$

Note that for $t \in [0, 1]$,

$$\bar{Z}^n(t) = \int_{[0,t] \times \Sigma'} \gamma \mu^n(d\theta, dx, d\gamma).$$

Since the trajectories of $\bar{Z}^n$ are Lipschitz continuous with a Lipschitz constant that is independent of $n$ and $\omega$, the limit $\bar{Z}$ is also Lipschitz continuous. Since $\mu(\{t\} \times \Sigma') = 0$ (a.e. in $t$), the weak convergence $\mu^n \to \mu$ and the last display imply

$$\bar{Z}(t) = \int_{[0,t] \times \Sigma'} \gamma \mu(d\theta, dx, d\gamma) \tag{21}$$

for all $t \in [0, 1]$. It follows from (18) that

$$\int_{[0,1]} \mathbf{1}\{\bar{X}^n(t) \in \mathcal{F}^{K,\emptyset}\} dR^n(t) = 0.$$

Thus $\mu^n(\Sigma_1) = 0$ for all $n \in \mathbb{N}$, and since $\Sigma_1$ is open relative to $\Sigma$ the weak convergence also implies

$$\mu(\Sigma_1) = 0. \tag{22}$$

Now by part 4 of Condition 4, $\gamma^n(t) \in d(\bar{X}^n(t)) \cup \{0\}$ for all $t \in [0, 1]$, and by (18) $(dR^n/dt)(t) = 0$ whenever $\gamma^n(t) = 0$. Thus

$$\int_{[0,1]} \mathbf{1}\{\gamma^n(t) \notin d(\bar{X}^n(t))\} dR^n(t) = 0,$$

which is the same as saying $\mu^n(\Sigma_2) = 0$. Since $\Sigma_2$ is open relative to $\Sigma$, the weak convergence also gives $\mu(\Sigma_2) = 0$.

To finish the proof of properties 4 and 5 of the SP, we use the uniform convergence $\bar{X}^n \to \bar{X}$, which for fixed $\Delta > 0$ implies that $\mu^n(\Sigma_3^\Delta) = 0$ for all sufficiently large $n$. Again using the weak convergence, $\mu(\Sigma_3^\Delta) = 0$ for all $\Delta > 0$. Sending $\Delta \to 0$, it follows that $\mu$ is supported on $\bar{\Sigma}(0, 1)$, where

$$\bar{\Sigma}(0, t) \doteq \bigcup_{\theta \in [0,t]} \bar{\Sigma}(\theta),$$

and

$$\bar{\Sigma}(t) \doteq \{(t, x, \gamma) \in \Sigma : x \notin \mathcal{F}^{K,\emptyset}, \gamma \in d(x), x = \bar{X}(t)\}.$$

From equation (21) we get

$$\bar{Z}(t) = \int_{\bar{\Sigma}(0,t)} \gamma \mu(d\theta, dx, d\gamma)$$

for $t \in [0, 1]$. Let $\lambda(t) = \mu([0, t] \times \Sigma')$. Then there exists a measurable mapping $t \in [0, 1] \mapsto \mu(dx, d\gamma|t) \in \mathcal{M}(\Sigma')$, such that for each Borel set $A$,

$$\mu([0, t] \times A) = \int_{[0,t]} \int_A \mu(dx, d\gamma|\theta) \lambda(d\theta).$$

It follows that
$$\bar{Z}(t) = \int_{[0,t]} \zeta(\theta)\lambda(d\theta),$$
where
$$\zeta(t) = \int_{\bar{\Sigma}(t)} \gamma \mu(dx, d\gamma|t).$$

Writing $\zeta(t) = a(t)\gamma(t)$, where $a(t) = \|\zeta(t)\|$ and $\gamma(t) = \zeta(t)/a(t)$ if $a(t) > 0$ and $\gamma(t) = 0$ otherwise, we have that $\gamma(t) \in d(X(t))$ $\lambda$-a.e. Hence

$$|\bar{Z}|(t) = \int_{[0,t]} a(\theta)\lambda(d\theta), \tag{23}$$

and consequently,

$$\bar{Z}(t) = \int_{[0,t]} \gamma(\theta)|\bar{Z}|(d\theta). \tag{24}$$

It follows from (23) that $|\bar{Z}|$ is absolutely continuous with respect to $\lambda$. Hence property 4 of the SP is implied by (22), and likewise property 5 is implied by (24).

# 5 Examples

In this section we verify all the conditions required by Theorem 4 for two interesting examples: a processor sharing model and the classical Jackson network. For both these systems the full model and all localized models will take the form assumed in Sections 2 to 4. For each example we must carry out the following tasks.

- Verify that the jump rates of the original model are constant in each facet.
- Check the communication condition (Condition 2).
- Verify the required uniqueness and perturbation properties of solutions to the system (8).
- Check that all corresponding SPs are regular whenever $c_v > 0$ for all $v \in V$.

## 5.1 The processor sharing model

Consider a queueing system that consists of a server and $N$ classes of customers that can be served, one at a time. There is one queue for each class, and the customers arrive and enter the queues according to their class. A column vector $f = (f_1, \ldots, f_N)^T$ satisfying $\sum_i f_i = 1$ is given, where the number $f_i > 0$ represents the minimal fraction of the overall service capacity guaranteed to class $i$ (see Figure 4). The probabilistic model is as follows. The arrivals are modeled as $N$ independent Poisson processes with rates $a_i > 0$, $i = 1, \ldots, N$. The service times of customers of class $i$ are exponential random variables with parameter $\sigma_i > 0$, independent of each other, of the service times of other classes, and of the arrival processes. When the server is free, service will be offered to a customer in one of the non-empty queues, say queue $i$, chosen at random, with probability proportional to $f_i\sigma_i$, independently of the other choices, service times and arrivals. If all queues are empty then no service takes place.

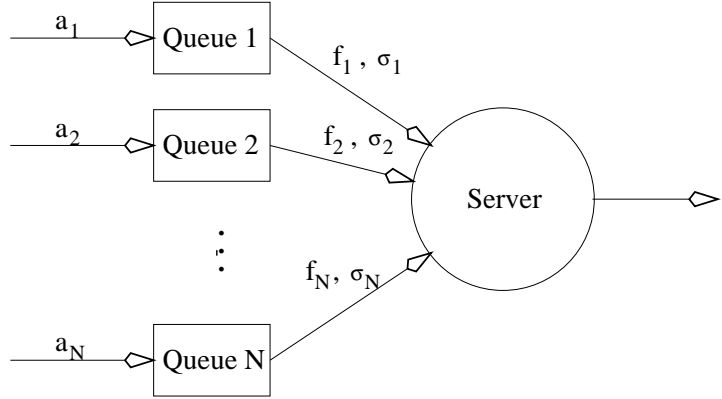

Figure 4: Processor sharing model

The vector defined by the number of customers of each class that are in the system, including those in the queues and the one being served, is a jump Markov process on $\mathbb{Z}_+^N$. Its jump intensity from $x \in \mathbb{Z}_+^N$ to $x + e_i$ is $a_i$ for $i = 1, \ldots, N$. Moreover, if $x \in \mathbb{Z}_+^N \cap \mathcal{F}^{\{1,\ldots,N\},I}$ and $I \neq \{1, ..., N\}$, then for $i \notin I$ the intensity of the jump from $x$ to $x - e_i$ is $\sigma_i f_i / f_{I^c}$, where $f_{I^c} = \sum_{i \in I^c} f_i$ and $I^c \doteq \{1,...,N\} \setminus I$. Thus the relative fraction of service offered to all nonempty queues is independent of the state of the system. Other models, such as models with modulated arrival rates and discrete time models could also be considered. However, the main point of the analysis is to show how one can deal with discontinuities in the statistical behavior in this multidimensional setting.

As noted in the last paragraph, the set $V$ of allowed jump directions is $\{\pm e_i, i = 1, \ldots, N\}$. In addition, the intensity function $r_{I,v}$ is given by $r_{I,e_i} = a_i$ for all facets $I$ and all $i = 1, \ldots, N$, and $r_{I,-e_i} = \sigma_i f_i / f_{I^c}$ if $i \notin I$ and $r_{I,-e_i} = 0$ otherwise. Thus Condition 1 holds.

We next consider the most involved condition, which is Condition 4 on the regularity of the associated SPs. This requires that we first identify the controlled processes that must be considered. Recall that these processes will have jump rates that are perturbed versions of the original jump rates. Given a solution of the system (8), the new jump rates will take the form $u_{I,v} = c_v r_{I,v}$. In particular, for any $v$ $u_{I,v}/r_{I,v}$ is independent of the facet $I$. We will show that Condition 4 holds for $K = \{1, \ldots, N\}$. As noted before Lemma 2, this implies that the condition also holds for all $K \subset \{1, \ldots, N\}$, i.e., all local models.

To simplify the notation we will write $c_i^{\pm} = c_{\pm e_i}$, $i = 1, \ldots, N$. The vectors that appear in Condition 4 can be calculated as follows. If $I \neq \{1, ..., N\}$, then

$$\sum_v c_v(r_{I,v} - r_{\emptyset,v})v = -\sum_{i \in I^c} c_i^- \sigma_i \frac{f_i}{f_{I^c}} e_i + \sum_{i=1}^N c_i^- \sigma_i f_i e_i.$$

For facets $\mathcal{F}^{\{1,\ldots,N\},I}$ of co-dimension 1, namely for $I = \{j\}$, it is easy to verify that the last display

reduces to
$$\sum_v c_v(r_{I,v} - r_{\emptyset,v})v = \frac{f_j}{1-f_j}C\Lambda(e_j - f),$$
where $C = \text{diag}(c_i^- : i = 1, \ldots, N)$ and $\Lambda = \text{diag}(\sigma_i : i = 1, \ldots, N)$. We define
$$d_j \doteq \frac{C\Lambda(e_j - f)}{\|C\Lambda(e_j - f)\|}$$
for $j = 1, \ldots, N$. If $n_j = e_j$, then obviously $\langle d_j, n_j \rangle > 0$, $j = 1, \ldots, N$. We next examine facets $\mathcal{F}^{\{1,\ldots,N\},I}$ of co-dimension higher than 1, (i.e., we let $I$ consist of more than 1 element), but still assume $I \neq \{1,\ldots,N\}$. In this case
$$\sum_v c_v(r_{I,v} - r_{\emptyset,v})v = \sum_{i \in I} \frac{f_i}{f_{I^c}} C\Lambda(e_i - f) = \sum_{i \in I} \alpha_i d_i,$$
where $\alpha_i = \|C\Lambda(e_i - f)\| f_i/f_{I^c}$. Hence if the expression in the latter display is non-zero then, when normalized, it belongs to $d(x)$ for any point $x$ in the facet $\mathcal{F}^{\{1,\ldots,N\},I}$.

Lastly, we must consider the case $I = \{1,\ldots,N\}$, which corresponds to the facet $\{0\}$. For all $i = 1, \ldots, N$ $r_{I,-e_i} = 0$, and hence one constraint direction at the origin should be
$$\sum_v c_v(r_{I,v} - r_{\emptyset,v})v = \sum_{i=1}^N c_i^- \sigma_i f_i e_i = C\Lambda f. \tag{25}$$

Note that the latter cannot be obtained by a linear combination of $d_i$, $i = 1,\ldots,N$, since for all $i$
$$\langle C\Lambda(e_i - f), C^{-1}\Lambda^{-1}(1,\ldots,1)^T \rangle = 0,$$
while
$$\langle C\Lambda f, C^{-1}\Lambda^{-1}(1,\ldots,1)^T \rangle = 1.$$
Following the treatment in [9], we supplement the set of normals $n_i = e_i$, $i = 1, \ldots, N$ with
$$n_{N+1} \doteq C\Lambda(1,\ldots,1)^T / \|C\Lambda(1,\ldots,1)^T\|,$$
thus introducing the extra constraint $\{x \in \mathbb{R}^N : \langle x, n_{N+1} \rangle \geq 0\}$. Note that the domain $\mathbb{R}_+^N$ is still given by $\cap_{i=1}^{N+1}\{x \in \mathbb{R}^N : \langle x, n_i \rangle \geq 0\}$. By also setting $d_{N+1} = n_{N+1}$, it is now possible to obtain the direction (25) as a linear combination of $d_i$, $i = 1, \ldots, N+1$ with nonnegative coefficients. We have therefore shown that parts 1, 2 and 4 of Condition 4 hold with $q = N + 1$. To complete the verification of Condition 4, all that remains is to check the regularity of the SP. Regularity for the case where $C\Lambda$ is the identity matrix is proved in [9]. Extending the proof to general $C\Lambda$ is easy since this SP can be obtained from the case of the identity matrix by a diagonal change of variable.

**Theorem 6 (Dupuis and Ramanan [9])** *Consider the SP associated with $\{(n_i, d_i), i = 1, \ldots, N+1\}$, where for $i = 1, \ldots, N$, $n_i = e_i$, and $d_i = (e_i - f)/\|e_i - f\|$, and $n_{N+1} = d_{N+1} = \sum_{i=1}^N e_i/\sqrt{N}$. Then the corresponding SM is regular.*

Now let $b_i > 0$, $i = 1, \ldots, N$ be fixed, and denote $B = \text{diag}(b_i : i = 1, \ldots, N)$. Then we have the following.

**Corollary 1** *Consider the SP associated with $\{(\tilde{n}_i, \tilde{d}_i), i = 1, \ldots, N+1\}$, where for $i = 1, \ldots, N+1$, $\tilde{n}_i = Bn_i/\|Bn_i\|$ and $\tilde{d}_i = Bd_i/\|Bd_i\|$, and $n_i$ and $d_i$ are as in Theorem 6. Then the corresponding SM is regular.*

**Proof:** Note first that for $i = 1, \ldots, N$, $\tilde{n}_i = e_i$, and $\tilde{n}_{N+1} = B(1, \ldots, 1)^T/\|B(1, \ldots, 1)^T\|$. Therefore the definition of $\tilde{n}_i$ is consistent with $G$, which is still given by $\cap_{i=1}^{N+1}\{x \in \mathbb{R}^N : \langle x, \tilde{n}_i \rangle \geq 0\}$. We show below that for any $\psi \in D_G([0,1] : \mathbb{R}^N)$ which is of bounded variation, $(\phi, \eta)$ solves the SP for $\psi$ with respect to $G$ and $d$ if and only if $(B\phi, B\eta)$ solves the SP for $B\psi$ with respect to $G$ and $\tilde{d}$. On the subset of $D_G([0,1] : \mathbb{R}^N)$ of paths of bounded variation, this would imply uniqueness and Lipschitz continuity of the SM associated with $G$ and $\tilde{d}$, based on these properties for the one associated with $G$ and $d$. The extension to $D_G([0,1] : \mathbb{R}^N)$ is as in [9].

Let us then assume that $(\phi, \eta)$ solves the SP for $\psi$ with respect to $G$ and $d$. It remains to check properties 1-5 of Definition 1 for $(B\phi, B\eta)$ and $B\psi$ with respect to $G$ and $\tilde{d}$. Properties 1-3 are immediate. Note that $\phi(\theta) \in \partial G$ if and only if $B\phi(\theta) \in \partial G$, and that $|\langle B\eta, e_i \rangle|$ is absolutely continuous with respect to $|\eta|$. Since we have that $\int_{[0,1]} \mathbf{1}\{\phi(\theta) \notin \partial G\} d|\eta|(\theta) = 0$, it follows that $\int_{[0,1]} \mathbf{1}\{B\phi(\theta) \notin \partial G\} d|\langle B\eta, e_i\rangle|(\theta) = 0$, $i = 1, \ldots, N$, and hence $\int_{[0,1]} \mathbf{1}\{B\phi(\theta) \notin \partial G\} d|B\eta|(\theta) = 0$. Thus property 4 follows.

Let $\gamma(\cdot)$ be as in property 5 of Definition 1, corresponding to $(\psi, \phi, \eta)$. Define $\tilde{\gamma}(t) = B\gamma(t)/\|B\gamma(t)\|$. Since $d(\phi(t)) = \tilde{d}(B\phi(t))$, $t \in [0,1]$, the definition of $\tilde{d}$ implies that $\tilde{\gamma}(t) \in \tilde{d}(B\phi(t))$ whenever $\gamma(t) \in d(\phi(t))$. Since $\eta(t) = \int_{[0,t]} \gamma(\theta) d|\eta|(\theta)$, we have

$$|B\eta|(t) = \int_{[0,t]} \|B\gamma(\theta)\| d|\eta|(\theta).$$

Thus $B\eta(t) = \int_{[0,t]} (B\gamma(\theta)/\|B\gamma(\theta)\|) d|B\eta|(\theta)$, and property 5 follows. The reverse direction is obtained similarly, using the transformation $B^{-1}$ instead of $B$. The proof of the corollary is therefore complete.
□

In view of the corollary, Condition 4 is verified, and we now turn Condition 2. It suffices to verify this condition for just the full model. We first identify for each $x \in \mathbb{Z}_+^N$ those directions $v$ for which $r(x, v) > 0$. Since $a_i > 0$ for $i = 1, \ldots, N$, we have that $r(x, e_i) > 0$, $i = 1, \ldots, N$. Recall that $I(x)$ is defined by $I(x) = \{i \in \{1, \ldots, N\} : \langle x, e_i \rangle = 0\}$. Since $\sigma_i f_i > 0$, $i = 1, \ldots, N$, we have that $r(x, -e_i) > 0$ for all $i \in \{1, \ldots, N\} \setminus I(x)$. It follows that

$$r(x, v) > 0 \text{ for all } x \in \mathbb{Z}_+^N \text{ and } v \in V \text{ for which } x + v \in \mathbb{Z}_+^N.$$

In other words, the set of possible jumps is always $V$, unless such a jump takes the process out of $\mathbb{Z}_+^N$. Suppose we are given points $x, y \in \mathbb{Z}_+^N$. Then it is easy to construct a sequence of points $\{x_0, x_1, \ldots, x_J\}$ with $x_0 = x$, $x_J = y$, $x_i - x_{i-1} \in V$ and $x_i \in \mathbb{Z}_+^N$ for all $i = 1, \ldots, J$, where $J = \sum_{k=1}^N |\langle y, e_k \rangle - \langle x, e_k \rangle|$. By the last display, we have that for any such sequence, $r(x_{i-1}, x_i - x_{i-1}) > 0$ for all $i = 1, \ldots, J$, and Condition 2 follows.

We next verify Condition 3. This condition must be verified for each local model, and so we fix $K \subset \{1,...,N\}$. Suppose the quadruple $(\{\rho_I\}_{I \subset K}, \{\bar{r}_v\}_{v \in V}, \{c_v\}_{v \in V}, \beta)$ satisfies (8), where $c_v > 0$, $v \in V$, and $\beta \in \mathcal{F}^{K,K}$. Since $r_{I,e_i} = a_i$, $I \subset K$, $i = 1, \ldots, N$, it follows from (8) that $\bar{r}_{e_i} = a_i$, $i = 1, \ldots, N$. Let us denote $\tau_i \doteq \bar{r}_{-e_i}$. Then for $i = 1, \ldots, N$ the last line of (8) implies

$$c_i^+ a_i - c_i^- \tau_i = \langle \beta, e_i \rangle, \tag{26}$$

and hence all $\{\bar{r}_v\}_{v \in V}$ are uniquely determined. Similarly, if $(\{\rho'_I\}_{I \subset K}, \{\bar{r}'_v\}_{v \in V}, \{c_v\}_{v \in V}, \beta')$ satisfies (8), and $\beta' \in \mathcal{F}^{K,K}$, then as before $\bar{r}'_{e_i} = a_i$, and the last display together with its primed analogue imply $c_i^- (\tau'_i - \tau_i) = \langle \beta - \beta', e_i \rangle$. Hence $|\bar{r}'_v - \bar{r}_v| \leq \|\beta - \beta'\| / \min_i c_i^-$, $v \in V$, and part 1 of Condition 3 follows.

As for part 2 of Condition 3, suppose that the quadruple $(\{\rho_I\}_{I \subset K}, \{\bar{r}_v\}_{v \in V}, \{c_v\}_{v \in V}, \beta)$ satisfies (8), where $\beta \in \mathcal{F}^{K,K}$, but we only know that $c_v \geq 0$, $v \in V$. If $\rho_\emptyset > 0$, we set $\rho'_I = \rho_I$, $I \subset K$. Otherwise, for some small $\kappa > 0$ we set $\rho'_\emptyset = \kappa$, and $\rho'_{I^*} = \rho_{I^*} - \kappa$ for some $I^*$ for which $\rho_{I^*} > \kappa$. For all $I \notin \{\emptyset, I^*\}$ we take $\rho'_I = \rho_I$. Since $r_{\emptyset,v} > 0$, $v \in V$, it follows that $\bar{r}'_v \doteq \sum_I \rho'_I r_{I,v} > 0$ for all $v \in V$. For $i = 1, \ldots, N$, we next define $c_i^{-\prime} = c_i^- + \kappa$ and

$$c_i^{+\prime} = \begin{cases} c_i^+ + \kappa & i \notin K, \\ c_i^{-\prime} \tau'_i / a_i & i \in K, \end{cases}$$

where as before $\tau'_i$ denotes $\bar{r}'_{-e_i} > 0$. With these definitions we have $c_i^{+\prime} a_i - c_i^{-\prime} \tau'_i = \langle \beta', e_i \rangle$, where $\beta' = \beta + \sum_{i \notin K} (\kappa a_i - \kappa \tau'_i) e_i \in \mathcal{F}^{K,K}$. Hence the quadruple $(\{\rho'_I\}_{I \subset K}, \{\bar{r}'_v\}_{v \in V}, \{c'_v\}_{v \in V}, \beta')$ satisfies (8). Since $\tau'_i > 0$ we have $c_i^{\pm\prime} > 0$. Finally, we observe that if $\kappa \to 0$ then $\rho'_I \to \rho_I$, $\bar{r}'_v \to \bar{r}_v$, $c'_v \to c_v$, and $\beta' \to \beta$. This proves that a slightly stronger statement than part 2 of Condition 3 holds, since $c'_v \to c_v$ for all $v \in V$, and in particular, $c'_v$ are bounded for $v \in \{w \in V : \bar{r}_w = 0\}$. This completes the verification of Condition 3. With all four conditions verified, the identification of the rate function for the processor sharing model is complete.

We have just shown that all the conditions of Theorem 4 hold for the model of this subsection. Let us rephrase the formula (17) in the current setting. The sum $\sum_v \bar{r}_v \ell(c_v)$ translates to $\sum_{i=1}^N [a_i \ell(c_i^+) + \tau_i \ell(c_i^-)]$. Also, the equation $\bar{r}_v = \sum_I \rho_I r_{I,v}$ in (8) translates as follows: $\bar{r}_{e_i} = a_i$, and

$$\tau_i = \bar{r}_{-e_i} = \sum_I \rho_I r_{I,-e_i} = \sum_{I \subset K : i \in I^c} \rho_I \frac{\sigma_i f_i}{f_{I^c}},$$

where $f_{I^c} = \sum_{j \in I^c} f_j$. The special form that Theorem 4 takes in this setting is therefore stated as follows.

**Theorem 7** *Let $a_i > 0$, $f_i > 0$ and $\sigma_i > 0$, $i = 1, \ldots, N$ be given, and consider the associated processor sharing model, as described in this subsection. Then the rate function for the local model corresponding to $K \subset \{1, \ldots, N\}$ takes the following form, for $\beta \in \mathcal{F}^{K,K}$:*

$$L(\beta) = \inf \left\{ \sum_{i=1}^N [a_i \ell(c_i^+) + \tau_i \ell(c_i^-)] : \begin{array}{l} c_i^+ a_i - c_i^- \tau_i = \langle \beta, e_i \rangle, \\ \tau_i = \sigma_i f_i \sum_{I \subset K : i \in I^c} \rho_I / f_{I^c}, \\ \rho_I \geq 0, \sum_I \rho_I = 1, c_i^\pm \geq 0 \end{array} \right\},$$

*where $f_{I^c} = \sum_{j \notin I} f_j$.*

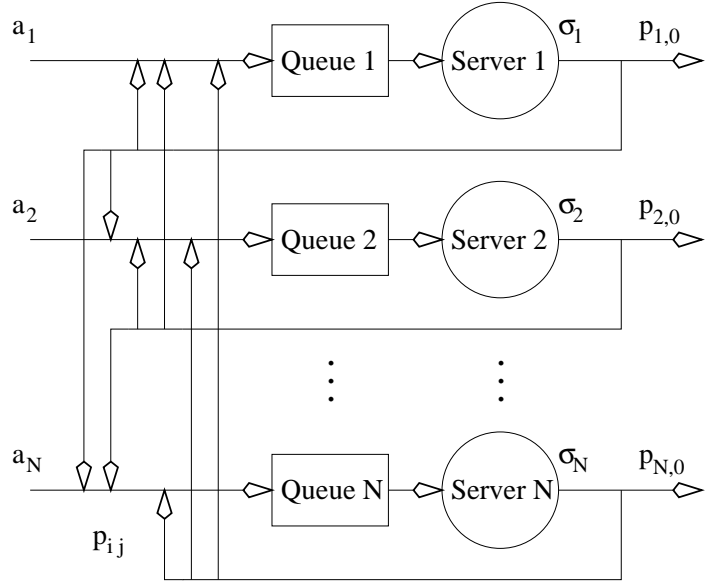

Figure 5: Jackson network

## 5.2 The Jackson network

Consider a queueing system in which customers of one class occupy $N$ nodes, where each node consists of a queue and a server. Customers arrive at a given queue from either other nodes or outside the system, and after being served, they may move to one of the $N$ queues (including the queue at the current node), or exit the system. The statistics of arrival, service time and routing variables depend on the node to which they correspond (see Figure 5). In particular, the arrivals are modeled as independent Poisson processes with rates $a_i \geq 0$, $i = 1, \ldots, N$, where $a_i > 0$ for at least one $i = 1, \ldots, N$. Service times are independent of each other and of the arrival processes, and for server $i$ they are exponential random variables with parameter $\sigma_i > 0$, $i = 1, \ldots, N$. The routing variables are independent of each other and of the arrivals and service times. We denote the probability that a customer leaving server $i$ is routed to queue $j$ by $p_{i,j}$, and the probability that they exit the system by $p_{i,0}$, for $i, j = 1, \ldots, N$. We assume that $p_{i,0} > 0$ for at least one $i = 1, \ldots, N$. Moreover, we assume that the sub-stochastic matrix $\{p_{i,j} : i, j = 1, \ldots, N\}$ is irreducible. Thus for every $i, j = 1, \ldots, N$ there exists a sequence of indices $i = i_0, i_1, \ldots, i_J = j$ such that $p_{i_{k-1}, i_k} > 0$ for $k = 1, \ldots, J$.

For this process, the vector defined by the number of customers at each node is a jump Markov process on $\mathbb{Z}_+^N$. The event of arrival of a customer to node $i$ corresponds to a jump in direction $e_i$, routing from node $i$ to node $j \neq i$ to a jump in direction $e_{i,j} \doteq e_j - e_i$, and exiting the system from node $i$ to a jump in direction $-e_i$. When a customer is routed from a node back to the same node, no jump occurs. The jump intensity from $x \in \mathbb{Z}_+^N$ in direction $e_i$ is therefore $a_i$, $i = 1, \ldots, N$. Moreover, if $x \in \mathbb{Z}_+^N \cap \mathcal{F}^{\{1,\ldots,N\},I}$ and $i \in I^c$ then the jump intensity in direction $e_{i,j}$ is $\sigma_i p_{i,j}$, while in direction $-e_i$ it is $\sigma_i p_{i,0}$. If $i \in I$ then both jumps have intensity zero. Thus Condition 1 holds.

Define
$$H^+ \doteq \{i : a_i > 0\}, \quad H^- \doteq \{i : p_{i,0} > 0\}, \quad H \doteq \{(i,j) : p_{i,j} > 0, i \neq j\}. \tag{27}$$

The set $V$ is given by $V = \{e_i, i \in H^+\} \cup \{-e_i, i \in H^-\} \cup \{e_{i,j}, (i,j) \in H\}$. Given $I \subset \{1,...,N\}$ the intensity function is given by $r_{I,e_i} = a_i$, $r_{I,e_{i,j}} = \sigma_i p_{i,j} \mathbf{1}\{i \notin I\}$, $r_{I,-e_i} = \sigma_i p_{i,0} \mathbf{1}\{i \notin I\}$. The jump intensities for the controlled processes we consider are given by $u_{I,v} = c_v r_{I,v}$, where $c_v > 0$, $v \in V$ are fixed.

As in Subsection 5.1, our most involved task is to verify Condition 4 for $K = \{1,\ldots,N\}$. We denote $c_i^\pm = c_{\pm e_i}$ and $c_{i,j} = c_{e_{i,j}}$, for $i,j = 1,\ldots,N$, $i \neq j$. For $I \subset \{1,...,N\}$ we have

$$\begin{aligned} \sum_v c_v(r_{I,v} - r_{\emptyset,v})v &= \sum_{i=1}^N \sum_{\substack{j=1\\j\neq i}}^N c_{i,j}(\sigma_i p_{i,j}\mathbf{1}\{i \notin I\} - \sigma_i p_{i,j})e_{i,j} - \sum_{i=1}^N c_i^-(\sigma_i p_{i,0}\mathbf{1}\{i \notin I\} - \sigma_i p_{i,0})e_i \\ &= -\sum_{\substack{i=1\\i\in I}}^N \sum_{\substack{j=1\\j\neq i}}^N c_{i,j}\sigma_i p_{i,j} e_{i,j} + \sum_{\substack{i=1\\i\in I}}^N c_i^-\sigma_i p_{i,0} e_i \end{aligned} \tag{28}$$

Therefore, we define the SP $\{(n_i, d_i), i = 1\ldots,q\}$ of Condition 4 by setting $q = N$, $n_i = e_i$, and

$$d_i = -\sum_{\substack{j=1\\j\neq i}}^N c_{i,j}\sigma_i p_{i,j} e_{i,j} + c_i^- \sigma_i p_{i,0} e_i \tag{29}$$

for $i = 1,\ldots,N$. Part 1 of Condition 4 obviously holds. By the irreducibility assumption, for every $i = 1,\ldots,N$ there exists $j \in \{0,1,...,N\}\setminus\{i\}$ for which $p_{i,j} > 0$. Hence

$$\langle d_i, n_i \rangle = \sum_{\substack{j=1\\j\neq i}}^N c_{i,j}\sigma_i p_{i,j} + c_i^- \sigma_i p_{i,0} > 0,$$

and part 2 of Condition 4 holds. To check part 4 of Condition 4, note that

$$\sum_v c_v(r_{I,v} - r_{\emptyset,v})v = \sum_{i\in I} d_i.$$

If the quantity in the last display is non-zero then it points in a direction in $d(x)$, where $x \in \mathcal{F}^{\{1,\ldots,N\},I}$.

To show regularity of the SP let us state the following result from [9], which is a slight generalization of a result in [10].

**Theorem 8** *Consider the SP associated with $G \subset \mathbb{R}^N$ and $\{(d_i, n_i), i = 1\ldots,N\}$, where the directions of constraint $\{d_i\}$ are linearly independent. Let the matrix $Q$ be defined by*

$$Q = [q_{i,j}] \doteq [|\delta_{i,j} - \langle d_i, n_j\rangle/\langle d_i, n_i\rangle|],$$

*where $\delta_{i,j} = 1$ if $i = j$ and equals zero otherwise. If $\sigma(Q) < 1$, then the corresponding SM is regular.*

We show that the directions of constraint $d_i$, $i = 1, \ldots, N$ defined by (29) are linearly independent. Suppose for some $\alpha_1, \ldots, \alpha_N$, one has $\sum_i \alpha_i d_i = 0$. Setting $c_{i,i} = 1$, $i = 1, \ldots, N$, we have by (29) that

$$\sum_{i=1}^{N} \alpha_i \left\{ c_i^- \sigma_i p_{i,0} + \sum_{j=1}^{N} c_{i,j} \sigma_i p_{i,j} \right\} e_i = \sum_{i,j=1}^{N} \alpha_i c_{i,j} \sigma_i p_{i,j} e_j. \tag{30}$$

Let

$$\gamma_i = c_i^- \sigma_i p_{i,0} + \sum_{j=1}^{N} c_{i,j} \sigma_i p_{i,j}, \quad i = 1, \ldots, N,$$

and note that $\gamma_i > 0$. Let also

$$\tilde{p}_{i,j} = c_{i,j} \sigma_i p_{i,j} / \gamma_i, \quad i, j = 1, \ldots, N,$$

and note that $[\tilde{p}_{i,j}]$ is a strictly sub-stochastic matrix, since $p_{i,0} > 0$ for some $i$, and hence 1 cannot be an eigenvalue. However, with this notation (30) becomes

$$\sum_{i=1}^{N} \alpha_i \gamma_i e_i = \sum_{i,j=1}^{N} \alpha_i \gamma_i \tilde{p}_{i,j} e_j,$$

which implies that $\alpha_i = 0$, $i = 1, \ldots, N$. Thus the $d_i$ are linearly independent. To calculate the $q_{i,j}$ of Theorem 8, note that for $i \neq j$, $\langle d_i, n_j \rangle = c_{i,j} \sigma_i p_{i,j}$, and thus $q_{i,i} = 0$, while

$$q_{i,j} = \frac{c_{i,j} \sigma_i p_{i,j}}{\sum_{k \neq i} c_{i,k} \sigma_i p_{i,k} + c_i^- \sigma_i p_{i,0}},$$

for $j \neq i$. From this it is evident that $Q = [q_{i,j}]$ is strictly sub-stochastic, and therefore $\sigma(Q) < 1$. Since Theorem 8 holds, so does Condition 4. As discussed prior to Lemma 2 it therefore also holds for all $K \subset \{1, \ldots, N\}$.

Before verifying Conditions 2 and 3 we fix $K \subset \{1, \ldots, N\}$ and consider the corresponding local process, which is a jump Markov process on $\mathcal{S} = \mathbb{Z}_+^{N,K}$. The intensity function is given by

$$r_{I,v} = \begin{cases} a_i & v = e_i, \\ \sigma_i p_{i,0} & v = -e_i, i \notin I, \\ \sigma_i p_{i,j} & v = e_{i,j}, i \notin I, \end{cases} \tag{31}$$

and zero otherwise. It is useful to notice that one also has the following expression for the intensity function: for $x \in \mathbb{Z}_+^{N,K}$, $v \in V$,

$$r(x,v) = r_{\emptyset,v} \mathbf{1}\{x + v \in \mathbb{Z}_+^{N,K}\}.$$

In other words, the process uses the generator of the original queueing model that applies in the facet $\emptyset$, except possibly at points on the "boundary" (i.e., a facet), in which case it simply deletes any jumps that would take it outside $\mathbb{Z}_+^{N,K}$.

We need verify Condition 2 only for the full model, and so take $K = \{1, \ldots, N\}$. This condition requires that we connect any two points in $\mathcal{S} = \mathbb{Z}_+^N$ (in terms of the intensity function being positive) with a linear bound on the length of the connecting sequence. We first show that it holds for any two points on any "simplex" of the form $\{x \in \mathbb{Z}_+^N : \sum_i \langle x, e_i \rangle = A\}$, and then extend to all of $\mathbb{Z}_+^N$. Within each simplex, we show the existence of a connecting sequence by first constructing one for each two "neighboring points" (i.e., points $x, y$ such that $y = x - e_i + e_j$), and then moving between any two points on the simplex along neighboring points.

Let $i, j \in \{1, \ldots, N\}$, and let $i \neq j$. Then by the irreducibility assumption there exists a sequence of indices $i = i_0, i_1, \ldots, i_J = j$, where $J \leq N$, for which $p_{i_{k-1}, i_k} > 0$ for $k = 1, \ldots, J$. Hence, for any $z \in \mathbb{Z}_+^N$ it follows from (31) that the sequence $z_k = z + e_{i_k}$ $k = 0, 1, \ldots, J$, is such that $r(z_k, z_{k+1} - z_k) > 0$ for $k = 1, \ldots, J$. Suppose that $x, y \in \mathbb{Z}_+^N$ are such that $y = x - e_i + e_j$. If we show that $z \doteq x - e_i \in \mathbb{Z}_+^N$, then it will follow that there exists a sequence $x = x_0, x_1, \ldots, x_J = y$, $x_k \in \mathbb{Z}_+^N$, $k = 1, \ldots, J$ such that $r(x_k, x_{k+1} - x_k) > 0$ for $k = 1, \ldots, J$. However, $\langle z, e_i \rangle \geq 0$ follows from $\langle z, e_i \rangle = \langle y, e_i \rangle$, and the fact that $y \in \mathbb{Z}_+^N$. By induction, if $x, y \in \mathbb{Z}_+^N$ are any points such that $\sum_i \langle x, e_i \rangle = \sum_i \langle y, e_i \rangle$, then there exists a connecting sequence whose length is at most $\|x - y\|_1 N$.

Now suppose that $x, y \in \mathbb{Z}_+^N$ and $\Delta \doteq \sum_i \langle y, e_i \rangle - \sum_i \langle x, e_i \rangle > 0$, and recall that $a_j > 0$ for some $j \in \{1, \ldots, N\}$. Hence one can consider the sequence $x, x + e_j, \ldots, x + \Delta e_j$ and use the fact that all the intensities $r(x, e_j), \ldots, r(x + (\Delta - 1)e_j, e_j)$ are positive. To construct a sequence from $x + \Delta e_j$ to $y$ along which the transition intensities are positive, one can use the argument of the last paragraph. For the case where $\Delta < 0$, recall that $p_{k,0} > 0$ for some $k \in \{1, \ldots, N\}$. Thus for all $z \in \mathbb{Z}_+^N$ one has $r(z + e_k, z) > 0$. Hence if $\Delta < 0$ one can first construct a sequence from $x$ to $y + |\Delta| e_k$ along which the transition intensities are positive, again using the argument of the last paragraph, and then move from $y + |\Delta| e_k$ to $y$. The total length of the sequence from $x$ to $y$ in all cases is at most $(N+1)\|x - y\|_1$, and Condition 2 therefore holds.

We begin the verification of Condition 3 by rewriting the last equation in (8) in accordance with (31). Let the set $K \subset \{1, \ldots, N\}$ be fixed. Then

$$\beta = \sum_{I,v} c_v \rho_I r_{I,v} v$$

$$= \sum_{i=1}^N c_i^+ a_i e_i - \sum_{I \subset K} \sum_{i \notin I} c_i^- \rho_I \sigma_i p_{i,0} e_i + \sum_{I \subset K} \sum_{i \notin I} \sum_{j \neq i} c_{i,j} \rho_I \sigma_i p_{i,j} (e_j - e_i).$$

For $i = 1, \ldots, N$ define

$$\tau_i \doteq \sum_{I : i \notin I} \rho_I. \tag{32}$$

Then

$$\beta - \sum_{i=1}^N c_i^+ a_i e_i = \sum_{i=1}^N \tau_i \left\{ -c_i^- \sigma_i p_{i,0} e_i + \sum_{j=1}^N c_{i,j} \sigma_i p_{i,j} (e_j - e_i) \right\} \tag{33}$$

$$= -\sum_{i=1}^N \tau_i d_i,$$

where the $d_i$ are as in (29). If we write $D$ for the matrix $(d_1, \ldots, d_N)$ and $\tau$ for the vector

$(\tau_1, \ldots, \tau_N)^T$, then the last display becomes

$$D\tau = -\beta + \sum_{i=1}^n c_i^+ a_i e_i. \tag{34}$$

It was proved immediately after the statement of Theorem 8 that $D$ is non-singular. It follows that $\tau$ depends continuously on $\beta$. Further, for $v = e_i$, $\bar{r}_v = a_i$, and for $v = -e_i$ or $v = e_{i,j}$,

$$\bar{r}_v = \tau_i r_{\emptyset,v}. \tag{35}$$

Hence $\bar{r}_v$ depends continuously on $\tau$, and part 1 of Condition 3 follows. Note that in addition to (34), $\tau$ is required to satisfy some further conditions, e.g., $\tau_i = 1$ for $i \in K^c$. Hence for a given $\beta$ and set $\{c_v\}_{v \in V}$ there may be no solution to (8). However, this does not affect the last argument, since Condition 3 only refers to properties of the solutions of (8) when they exist.

Verifying part 2 of Condition 3, which is straightforward but detailed, is postponed to the appendix.

Having verified all assumptions of Theorem 4, we now rephrase it is the context of the Jackson network example. Let us return to (34) and recall that

$$D = [d_{i,j}], \quad d_{i,j} = \langle d_i, e_j \rangle = \begin{cases} \sum_{\substack{k=1 \\ k \neq i}}^N c_{i,k} \sigma_i p_{i,k} + c_i^- \sigma_i p_{i,0} & i = j \\ -c_{i,j} \sigma_i p_{i,j} & i \neq j \end{cases} \tag{36}$$

is a non-singular matrix if $c_i^-, c_{i,j} > 0$. In addition, we let

$$C^+ = \mathrm{diag}(c_i^+ : i = 1, \ldots, N), \quad a = (a_1, \ldots, a_N)^T. \tag{37}$$

Note that (8) corresponds to the following set of conditions:

$$\begin{cases} \rho_I \geq 0, \sum_I \rho_I = 1, \\ \tau_i = \sum_{I:i \notin I} \rho_I, \\ c_i^\pm, c_{i,j} \geq 0, \\ D\tau = C^+ a - \beta. \end{cases}$$

We then have the following.

**Theorem 9** *Let $a_i \geq 0$, $\sigma_i > 0$, $p_{i,j} \geq 0$, $i = 1, \ldots, N$, $j = 0, 1, \ldots, N$ be given, such that $\sum_{j=0}^N p_{i,j} = 1$, $i = 1, \ldots, N$, $[p_{i,j}]_{i,j \in \{1,\ldots,N\}}$ is irreducible, and $a_i > 0$, $p_{j,0} > 0$ for some $i, j \in \{1, \ldots, N\}$. Consider the associated Jackson network, as described above. Then the rate function for the local model, corresponding to $K \subset \{1, \ldots, N\}$, is as follows for $\beta \in \mathcal{F}^{K,K}$:*

$$L(\beta) = \inf \left\{ \sum_{i=1}^N \left[ a_i \ell(c_i^+) + \tau_i \sigma_i p_{i,0} \ell(c_i^-) + \tau_i \sigma_i \sum_{\substack{j=1 \\ j \neq i}}^N p_{i,j} \ell(c_{i,j}) \right] : \begin{array}{l} \rho_I \geq 0, \sum_I \rho_I = 1, \\ \tau_i = \sum_{I:i \notin I} \rho_I, \\ c_i^\pm, c_{i,j} \geq 0, \\ D\tau = C^+ a - \beta \end{array} \right\},$$

*where $D$, $C^+$ and $a$ are as in (36) and (37).*

# 6  Appendix

**Proof of part 2 of Theorem 2:**

The inequality $q_M^n(y;\beta,\epsilon) \leq q^n(y;\beta,\epsilon)$ holds for any $\beta \in \mathbb{R}^N$, $\epsilon > 0$ and $y \in \mathcal{S}^n$. Therefore

$$\lim_{M\to\infty}\lim_{\epsilon\to 0}\lim_{\delta\to 0}\limsup_{n\to\infty}\sup_{\{y\in\mathcal{S}^n:\|y\|\leq\delta\}} q_M^n(y;\beta,\epsilon) \leq \lim_{\epsilon\to 0}\lim_{\delta\to 0}\limsup_{n\to\infty}\sup_{\{y\in\mathcal{S}^n:\|y\|\leq\delta\}} q^n(y;\beta,\epsilon).$$

We show below that for all sufficiently small $\epsilon > 0$ there is $\eta = \eta(\epsilon) > 0$, $n_0 = n_0(\epsilon) < \infty$, and $M_0$ which depends only on $\beta$ such that

$$q^n(y;\beta,\epsilon) \leq q_M^n(y;\beta,\eta) \tag{38}$$

for all $M > M_0$ and $n \geq n_0$. Thus part 2 of Theorem 2 will follow from part 1 of the same theorem.

Consider the stopping time

$$T^n = \inf\{t \in [0,1] : \|\Xi^n(t) - t\beta\| \geq \eta\}$$

(where $\inf \emptyset = \infty$), and let $\mathcal{F}_{T^n}$ denote the stopping $\sigma$-field associated with $T^n$. It follows from (4) that

$$q_M^n(y;\beta,\eta) \geq \inf \bar{E}_y^n\left\{\int_0^{T^n\wedge 1}\sum_v r(\Xi^n(t),v)\ell\left[\frac{u^n(\Xi^n(t),v,t)}{r(\Xi^n(t),v)}\right]dt + M\,\mathbf{1}\{T^n<\infty\}\right\}, \tag{39}$$

where the infimum is over all admissible controls $u^n$ and associated controlled processes $\Xi^n$. The inequality is due to the fact that the upper limit of integration 1 has been replaced by $T^n \wedge 1$. Also, according to Theorem 1

$$q^n(y;\beta,\epsilon) = \inf \bar{E}_y^n\left\{\int_0^1\sum_v r(\Xi^n(t),v)\ell\left[\frac{u^n(\Xi^n(t),v,t)}{r(\Xi^n(t),v)}\right]dt + g(\Xi^n;\beta,\epsilon)\right\},$$

where the infimum is over all admissible controls $u^n$ and associated controlled processes $\Xi^n$. By part 1 of Theorem 2 for sufficiently small $\epsilon > 0$ there are $\delta = \delta(\epsilon) > 0$ and $n_0 = n_0(\epsilon) < \infty$ such that for all $n > n_0$ and all $y \in \mathcal{S}^n$, $\|y\| \leq \delta$ there exists an admissible control $\hat{u}^n$ for which

$$\bar{E}_y^n\left\{\int_0^1\sum_v r(\Xi^n(t),v)\ell\left[\frac{\hat{u}^n(\Xi^n(t),v,t)}{r(\Xi^n(t),v)}\right]dt + g(\Xi^n;\beta,\epsilon)\right\} \leq L(\beta) + 1 \doteq M_0 < \infty.$$

We will need to define an analogous control problem for processes whose value is specified at times $s \in [0,1]$ (rather than just $s = 0$). Accordingly, let $\bar{E}_{y,s}^n$ denote conditioning on $\Xi^n(s) = y$, and define $g_s(\phi;\beta,\epsilon)$ to be $\infty$ if $\|\phi(t) - t\beta\| \geq \epsilon$ for some $t \in [s,1]$, and zero otherwise. Then for sufficiently small $\epsilon > 0$ there are $\delta = \delta(\epsilon) > 0$ and $n_0 = n_0(\epsilon) < \infty$ such that for all $s \in [0,1]$, all $n > n_0$, and all $y \in \mathcal{S}^n$, $\|y - s\beta\| \leq \delta$, there exists an admissible control $\hat{u}^n(\cdot,\cdot,\cdot;s,y)$ for which

$$\bar{E}_{y,s}^n\left\{\int_s^1\sum_v r(\Xi^n(t),v)\ell\left[\frac{\hat{u}^n(\Xi^n(t),v,t;s,y)}{r(\Xi^n(t),v)}\right]dt + g_s(\Xi^n;\beta,\epsilon)\right\} \leq (1-s)L(\beta) + 1 \leq M_0.$$

Owing to the continuous dependence of the distribution of $\Xi^n$ on the control, $\hat{u}^n$ can be selected so that it is measurable in all variables.

We now define a composite control that will be used in the variational representation for $q^n(y; \beta, \epsilon)$. Let $\eta \in (0, \epsilon \wedge \delta)$. If $u^n$ is an arbitrary admissible control and $\hat{u}^n$ is as above, then let
$$\tilde{u}^n(v, t) = \begin{cases} u^n(\Xi^n(t), v, t) & \text{for } t \in [0, T^n \wedge 1), \\ \hat{u}^n(\Xi^n(t), v, t; T^n, \Xi^n(T^n)) & \text{for } t \in [T^n \wedge 1, 1], \end{cases}$$
This control is exactly the same as $u^n$ up until the first time that the controlled process leaves the $\eta$-neighborhood of $\beta t$. At that time the control $\hat{u}^n$ takes over, which will keep the process within the $\epsilon$-neighborhood of $t\beta$ for the remaining time with a total cost (over the remaining time) of less than $M_0$. Now since the controlled process stays within the $\epsilon$-neighborhood of $t\beta$ for all $t \in [0,1]$, the variational characterization of $q^n(y, \beta, \epsilon)$ implies

$$q^n(y; \beta, \epsilon) \leq \inf \left\{ \bar{E}_y^n \int_0^{T^n \wedge 1} \sum_v r(\Xi^n(t), v) \ell \left[ \frac{u^n(\Xi^n(t), v, t)}{r(\Xi^n(t), v)} \right] dt \right.$$
$$\left. + \bar{E}_y^n \bar{E}_y^n \left[ \int_{T^n \wedge 1}^1 \sum_v r(\Xi^n(t), v) \ell \left[ \frac{\hat{u}^n(\Xi^n(t), v, t; T^n, \Xi^n(T^n))}{r(\Xi^n(t), v)} \right] dt \Big| \mathcal{F}_{T^n} \right] \right\}.$$

Strictly speaking, the control $\bar{u}^n$ is not of feedback form, but instead belongs to the larger class of nonanticipating controls. However, using the fact that the infimum over these two classes is the same the inequality holds as stated. We next observe that the conditional expectation in the last display is bounded above by $M_0 \bar{P}_y^n \{T^n < \infty\}$. It follows from (39) that for any given admissible control $u^n$ the right hand side is a lower bound for the corresponding cost in the variational representation for $q_M^n(y; \beta, \eta)$, at least if $M \geq M_0$. Since the last display holds for all admissible controls $u^n$, (38) follows.
$\square$

**Proof of part 2 of Condition 3 for the Jackson network:**

We now verify part 2 of Condition 3 in the setting of Subsection 5.2. We first prove the following.

**Lemma 6** *Let $K \subset \{1, \ldots, N\}$ be fixed. Then given any numbers $0 \leq \tau_i \leq 1$, $i \in K$, one can associate numbers $\rho_I \geq 0$, $I \subset K$, $\sum_{I \subset K} \rho_I = 1$, such that (32) holds for $i \in K$, and such that if $i \notin K$ then (32) implies $\tau_i = 1$.*

**Proof:** Let $0 \leq \tau_i \leq 1$, $i \in K$ be given. For $I \subset K$, define
$$\rho_I = \left[ \prod_{i \in K \setminus I} \tau_i \right] \left[ \prod_{i \in I} (1 - \tau_i) \right].$$
Since for any finite set of numbers $\{\alpha_i\}_{i \in J}$ one has $\sum_{J' \subset J} [\prod_{i \in J \setminus J'} \alpha_i][\prod_{i \in J'} (1 - \alpha_i)] = 1$, (32) holds, as well as $\sum_{I \subset K} \rho_I = 1$. If $i \notin K$, then by the last sentence the right hand side of (32) is one.
$\square$

Suppose we are given a quadruple $(\{\rho_I\}_{I \subset K}, \{\bar{r}_v\}_{v \in V}, \{c_v\}_{v \in V}, \beta)$ satisfying (8), with $\beta \in \mathcal{F}^{K,K}$. We let $\tau_i$, $i \in K$ be defined by (32) and observe that $\bar{r}_v$, $v \in V$ must then be given by (35).

Returning to the notation of (27), for any $i \in H^-$ ($i \in H^+$, respectively) for which $c_i^- = 0$ ($c_i^+ = 0$), we fix a sequence of indices $j = k_0^{-,i}, k_1^{-,i}, \ldots, k_{J_i^-}^{-,i} = i$ ($i = k_0^{+,i}, k_1^{+,i}, \ldots, k_{J_i^+}^{+,i} = j$), where $j \in H^+$ ($j \in H^-$), and $(k_{\ell-1}^i, k_\ell^i) \in H$ for $\ell = 1, \ldots, J_i^-$ ($\ell = 1, \ldots, J_i^+$). It is possible that a sequence consists of just one element, in which case the latter requirement is irrelevant. Similarly, for any $(i,j) \in H$ for which $c_{i,j} = 0$, we fix a sequence $j = k_0^{i,j}, k_1^{i,j}, \ldots, k_{J_{i,j}}^{i,j} = i$ along which we also have $(k_{\ell-1}^{i,j}, k_\ell^{i,j}) \in H$, $\ell = 1, \ldots, J_{i,j}$. We assume, without loss of generality, that each of the above sequences is self-avoiding, namely consists of distinct indices.

In order to define $c_i^{\pm\prime}$, $c_{i,j}{}'$ and $\tau_i'$, $i, j \in \{1, \ldots, N\}$, $i \neq j$ we first introduce some notation. Let

$$W_{i,j} \doteq \#\{(p,q) : (i,j) = (k_{m-1}^{p,q}, k_m^{p,q}), \text{ some } m = 1, \ldots, J_{p,q}\},$$

$$W_{i,j}^+ \doteq \#\{p : (i,j) = (k_{m-1}^{+,p}, k_m^{+,p}), \text{ some } m = 1, \ldots, J_p^+\},$$

$$W_{i,j}^- \doteq \#\{p : (i,j) = (k_{m-1}^{-,p}, k_m^{-,p}), \text{ some } m = 1, \ldots, J_p^-\},$$

$$W_i^+ \doteq \#\{p : i = k_{J_p^+}^{+,p}\},$$

$$W_i^- \doteq \#\{p : i = k_0^{-,p}\}.$$

The following are solely consequences of the definitions in the last two paragraphs. The first equality in each line is due to the fact that summation along all paths can be performed by first grouping the contributions of all paths to each index $(i,j)$ and $i$, and then summing on these indices. The second inequality in each line is due to the cancellation of all the middle terms along each path.

$$\sum_{i,j} W_{i,j}(e_j - e_i) = \sum_{\substack{(p,q) \in H: \\ c_{p,q}=0}} \sum_{m=1}^{J_{p,q}} (e_{k_m^{p,q}} - e_{k_{m-1}^{p,q}}) = \sum_{\substack{(p,q) \in H: \\ c_{p,q}=0}} e_p - e_q,$$

$$\sum_{i,j} W_{i,j}^-(e_j - e_i) + \sum_i W_i^- e_i = \sum_{\substack{p \in H^-: \\ c_p^-=0}} \left\{ \sum_{m=1}^{J_p^-} (e_{k_m^{-,p}} - e_{k_{m-1}^{-,p}}) + e_{k_0^{-,p}} \right\} = \sum_{\substack{p \in H^-: \\ c_p^-=0}} e_p,$$

$$\sum_{i,j} W_{i,j}^+(e_j - e_i) - \sum_i W_i^+ e_i = \sum_{\substack{p \in H^+: \\ c_p^+=0}} \left\{ \sum_{m=1}^{J_p^+} (e_{k_m^{+,p}} - e_{k_{m-1}^{+,p}}) - e_{k_{J_p^+}^{+,p}} \right\} = \sum_{\substack{p \in H^+: \\ c_p^+=0}} -e_p.$$

Let now $\kappa > 0$ be a number. If $i$ is such that $\tau_i > 0$, we let $\tau_i' = \tau_i$,

$$c_{i,j}' = c_{i,j} + \kappa \frac{\mathbf{1}\{c_{i,j} = 0\} + W_{i,j}^+ + W_{i,j}^- + W_{i,j}}{\tau_i \sigma_i p_{i,j}},$$

and

$$c_i^{-\prime} = c_i^- + \kappa \frac{\mathbf{1}\{c_i^- = 0\} + W_i^+}{\tau_i \sigma_i p_{i,0}}.$$

If $i$ is such that $\tau_i = 0$, we let $\tau_i' = \kappa$,

$$c_{i,j}' = \frac{\mathbf{1}\{c_{i,j} = 0\} + W_{i,j} + W_{i,j}^+ + W_{i,j}^-}{\sigma_i p_{i,j}},$$

and

$$c_i^- = \frac{\mathbf{1}\{c_i^- = 0\} + W_i^+}{\sigma_i p_{i,0}}.$$

Also, for $i \in H^+$ we let

$$c_i^{+\prime} = c_i^+ + \kappa \frac{\mathbf{1}\{c_i^+ = 0\} + W_i^-}{a_i}.$$

We now show that with the above definitions, (33) still holds true if one replaces $(\tau_i, c_i^\pm, c_{i,j})$ by $(\tau_i', c_i^{\pm\prime}, c_{i,j}')$, with the same $\beta$. In fact,

$$\sum_{i=1}^N c_i^{+\prime} a_i e_i + \sum_{i=1}^N \tau_i' \left\{ -c_i^{-\prime} \sigma_i p_{i,0} e_i + \sum_{j=1}^N c_{i,j}' \sigma_i p_{i,j}(e_j - e_i) \right\}$$
$$- \sum_{i=1}^N c_i^+ a_i e_i - \sum_{i=1}^N \tau_i \left\{ -c_i^- \sigma_i p_{i,0} e_i + \sum_{j=1}^N c_{i,j} \sigma_i p_{i,j}(e_j - e_i) \right\}$$
$$= \sum_i [\mathbf{1}\{c_i^+ = 0\} + W_i^-] \kappa e_i$$
$$\quad - \sum_i [\mathbf{1}\{c_i^- = 0\} + W_i^+] \kappa e_i$$
$$\quad + \sum_{i,j} [\mathbf{1}\{c_{i,j} = 0\} + W_{i,j} + W_{i,j}^+ + W_{i,j}^-] \kappa (e_j - e_i)$$
$$= \sum_{\substack{p \in H^+:\\ c_p^+ = 0}} -\kappa e_p + \sum_{\substack{p \in H^-:\\ c_p^- = 0}} \kappa e_p + \sum_{\substack{(p,q) \in H:\\ c_{p,q} = 0}} \kappa(e_p - e_q)$$
$$\quad + \sum_i \kappa \mathbf{1}\{c_i^+ = 0\} e_i - \sum_i \kappa \mathbf{1}\{c_i^- = 0\} e_i + \sum_{i,j} \kappa \mathbf{1}\{c_{i,j} = 0\}(e_j - e_i)$$
$$= 0.$$

Hence (33) holds for $(\tau_i', c_i^{\pm\prime}, c_{i,j}')$. Based on $\tau_i'$, $i = 1, \ldots, N$ one can now define $\rho_I'$, $I \subset K$ as in Lemma 6, and $\bar{r}_v'$ as in (35), and conclude that the quadruple $(\{\rho_I'\}_{I \subset K}, \{\bar{r}_v'\}_{v \in V}, \{c_v'\}_{v \in V}, \beta)$ satisfies (8). By definition we have that for all $\kappa > 0$, $c_v > 0$, $v \in V$, and also that as $\kappa \to 0$, $\tau_i' \to \tau_i$, $i = 1, \ldots, N$, and $c_v' \to c_v$ for $v \in \{w \in V : \bar{r}_v \neq 0\}$ while $c_v'$ remains bounded for $v \in \{w \in V : \bar{r}_v = 0\}$. Consequently, part 2 of Condition 3 holds.
□

**Proof of Lemma 2:** For $x = 0$ there is nothing to prove, hence the contrary is assumed throughout the proof. By assumption there is $K_1 < \infty$ such that for all bounded variation $\psi_1, \psi_2 \in D_G([0, T] : \mathbb{R}^N)$

$$\sup_{t \in [0,T]} \|\Gamma(\psi_1)(t) - \Gamma(\psi_2)(t)\| \leq K_1 \sup_{t \in [0,T]} \|\psi_1(t) - \psi_2(t)\|. \qquad (40)$$

Let
$$\tilde{G} = \bigcap_{i \in I(x)} \{y \in \mathbb{R}^N : \langle y, n_i \rangle \geq 0\}.$$

We first show that existence of solutions to the SP on $\tilde{G}$ follows from existence of solutions to the SP on $G$. Let $\tilde{\psi} \in D_{\tilde{G}}([0,T] : \mathbb{R}^N)$ be of bounded variation, and define

$$a \doteq \max_{i \notin I(x)} \sup_{t \in [0,T]} \frac{K_1 \|\tilde{\psi}(t) - \tilde{\psi}(0)\| + \|\tilde{\psi}(0)\| + 1}{\langle x, n_i \rangle}.$$

Let $\psi(\cdot) = \tilde{\psi}(\cdot) + ax$, and denote by $(\phi, \eta)$ the solution to the SP for $\psi$ with respect to $G$ and $d_i$, $i = 1, \ldots, q$. Let also $\tilde{\phi}(\cdot) = \phi(\cdot) - ax$ and $\tilde{\eta} = \eta$. The definition of $a$ implies that $\tilde{\phi}(t) \in \tilde{G}$ for $t \in [0,T]$. Moreover, for $i \in I(x)$ and $t \in [0,T]$, $\langle \tilde{\phi}(t), n_i \rangle = \langle \phi(t), n_i \rangle$. We also show below that for $i \notin I(x)$, $t \in [0,T]$,

$$\langle \phi(t), n_i \rangle > 0. \tag{41}$$

It follows from the last two assertions that for all $t \in [0,T]$ and $i \in I(x)$, $\langle \phi(t), n_i \rangle = 0$ if and only if $\langle \tilde{\phi}(t), n_i \rangle = 0$, and that for all $t \in [0,T]$, $\phi(t) \in \partial G$ if and only if $\tilde{\phi}(t) \in \partial \tilde{G}$. Hence $(\tilde{\phi}, \tilde{\eta})$ solve the SP for $\tilde{\psi}$ with respect to $\tilde{G}$ and $d_i$, $i \in I(x)$.

To show that (41) holds for $i \notin I(x)$, use (40) and the definition of $a$ to conclude

$$\begin{aligned}
\langle \phi(t), n_i \rangle &\geq \langle \phi(0), n_i \rangle - \|\phi(t) - \phi(0)\| \\
&\geq \langle \phi(0), n_i \rangle - K_1 \|\psi(t) - \psi(0)\| \\
&= \langle \tilde{\psi}(0), n_i \rangle + a \langle x, n_i \rangle - K_1 \|\tilde{\psi}(t) - \tilde{\psi}(0)\| \\
&\geq 1.
\end{aligned}$$

This completes the proof of existence of solutions to the SP on $\tilde{G}$.

We denote by $\Gamma$ (respectively, $\tilde{\Gamma}$) the corresponding SM on $G$ (respectively, $\tilde{G}$), and prove that Lipschitz continuity of $\tilde{\Gamma}$ follows from (40). Let $\tilde{\psi}_1, \tilde{\psi}_2 \in D_{\tilde{G}}([0,T] : \mathbb{R}^N)$ be of bounded variation, and for $j = 1, 2$ let $(\tilde{\phi}_j, \tilde{\eta}_j)$ solve the SP for $\tilde{\psi}_j$ with respect to $\tilde{G}$ and $d_i, i \in I(x)$. Define

$$\bar{a} \doteq \max_{j=1,2} \max_{i \notin I(x)} \sup_{t \in [0,T]} \frac{\|\tilde{\phi}_j(t)\| + 1}{\langle x, n_i \rangle}.$$

For $j = 1, 2$ let $\psi_j(\cdot) = \tilde{\psi}_j(\cdot) + \bar{a}x$, $\phi_j(\cdot) = \tilde{\phi}_j(\cdot) + \bar{a}x$, and $\eta_j = \tilde{\eta}_j$. For $j = 1, 2$, $i \notin I(x)$ and $t \in [0,T]$ we have

$$\begin{aligned}
\langle \phi_j(t), n_i \rangle &\geq \langle \tilde{\phi}_j(t), n_i \rangle + \sup_{\theta \in [0,T]} \|\tilde{\phi}_j(\theta)\| + 1 \\
&\geq 1.
\end{aligned}$$

It follows that $\phi_j(t) \in G$ for $j = 1, 2$ and $t \in [0,T]$, and moreover, that for all $i \in I(x)$, $\langle \tilde{\phi}_j(t), n_i \rangle = 0$ if and only if $\langle \phi_j(t), n_i \rangle = 0$, and that $\phi_j(t) \in \partial G$ if and only if $\tilde{\phi}_j(t) \in \partial \tilde{G}$. Because of this, $(\phi_j, \eta_j)$ solves the SP for $\psi_j$ with respect to $G$ and $d_i, i = 1, \ldots, q$, $j = 1, 2$. Hence

$$\sup_{t \in [0,T]} \|\tilde{\phi}_1(t) - \tilde{\phi}_2(t)\| \leq K_1 \sup_{t \in [0,T]} \|\tilde{\psi}_1(t) - \tilde{\psi}_2(t)\|,$$

and regularity of $\tilde{\Gamma}$ follows.

□

1991 *Mathematics Subject Classification.* Primary 60F10, 60K25; Secondary 93E20, 60F17.